\documentclass[12pt]{article}
\usepackage{kotex}
\usepackage{amsmath,amssymb,amsthm, mathrsfs, amsfonts}
\usepackage{latexsym}
\usepackage{hyperref}
\usepackage[bf, small]{titlesec}
\usepackage{authblk} 
\usepackage{booktabs}
\usepackage{array}

\theoremstyle{plain}
\newtheorem{Thm}{Theorem}[section]
\newtheorem{Lem}[Thm]{Lemma}
\newtheorem{Prop}[Thm]{Proposition}

\theoremstyle{definition}

\newtheorem{Ex}[Thm]{Example}

\usepackage{standalone}
\usepackage{tikz}
\usetikzlibrary{arrows,positioning,matrix,fit,backgrounds,shapes,shapes.geometric, intersections,decorations.pathreplacing,calc}

\usetikzlibrary{shapes.callouts,decorations.pathmorphing} 
\usetikzlibrary{shadows} 
\usepackage{calc} 
\usetikzlibrary{decorations.markings} 
\tikzstyle{vertex}=[circle, draw, inner sep=0pt, minimum size=6pt] 

\usetikzlibrary{arrows.meta}

\usepackage{mathtools}

\newcounter{statement}
\newcommand{\statement}[2]{%
\begin{equation}\refstepcounter{statement}\tag{S\thestatement}\label{#1}%
\parbox{\dimexpr\linewidth-4em}{#2}%
\end{equation}%
}

\title{On the convergence of the graph sequence $\left\{ C^m(D) \right\}_{m=1}^{\infty}$ for a multipartite tournament $D$}

\author[1]{\small Ji-Hwan Jung}
\author[2]{\small Suh-Ryung Kim}
\author[2]{\small Hyesun Yoon}
\affil[1]{\footnotesize Department of Mathematics Education, Chinju National University of Education, Jinju 52673}
\affil[2]{\footnotesize Department of Mathematics Education, Seoul National University, Seoul 08826}
\affil[ ]{\footnotesize\textit{jhjung@cue.ac.kr, srkim@snu.ac.kr, magisakura@snu.ac.kr}}
\date{}
\begin{document}
\maketitle
\begin{abstract}
Given a positive integer $m$, the {\em $m$-step competition graph} of a digraph $D$, denoted by $C^m(D)$, has the same vertex set as $D$ and has an edge between vertices $u$ and $v$ if and only if there exists a vertex $w$ such that there exist directed walks of length $m$ from $u$ to $w$ and from $v$ to $w$, respectively.
In this paper, we completely characterize the convergence of $\{C^m(D)\}_{m=1}^{\infty}$ for a multipartite tournament $D$ based on the last nontrivial strong component of $D$.
Furthermore, not only do we determine the limit in the case of convergence, but also in the event of divergence, we specify how $C^m(D)$ changes periodically depending on the value of $m$.
Our results extend the work of Jung~{\em et al.} [On the limit of the sequence $\{C^m (D)\}_{m=1}^{\infty}$ for a
  multipartite tournament $D$. {\em Discrete Appl. Math.}, 340:1--13, 2023] which addresses the case of the last strong component being nontrivial, thereby completing the convergence analysis of  $\{C^m(D)\}_{m=1}^{\infty}$ for a multipartite tournament $D$.
  Our results can also be expressed in terms of matrix sequence $\{A^m(A^T)^m\}_{m=1}^{\infty}$ for the adjacency matrix $A$ of $D$ and this part is also covered in the text.
\end{abstract}


\section{Introduction}
The underlying graph of each digraph in this paper is assumed to be simple unless otherwise mentioned.

Given vertex sets $X$ and $Y$ of a digraph, we use the notation $X \to Y$ for ``there is an arc $(x,y)$ for each $x \in X$ and for each $y \in Y$'' in the digraph.
For simplicity, we omit the braces if $X$ or $Y$ is a singleton set.

Let $D$ be a digraph and $m$ be a positive integer.
A vertex $y$ is an {\em $m$-step prey} of a vertex $x$ if and only if there exists a directed walk from $x$ to $y$ of length $m$.
The {\em $m$-step competition graph} of a digraph $D$, denoted by $C^m(D)$, has the same vertex set as $D$ and has an edge between vertices $u$ and $v$ if and only if there exists an $m$-step common prey of $u$ and $v$ in $D$.
The notion of an $m$-step competition graph introduced by Cho~{\em et al.}~\cite{cho2000m} is a generalization of the competition graph introduced by Cohen~\cite{cohen1968interval} (the {\em competition graph} of a digraph $D$ is $C^1(D)$).
Since its introduction, an $m$-step competition graph  has been extensively studied (see, for example, \cite{belmont2011complete,choi2023digraphs,helleloid2005connected,ho2005m,park2011m,zhao2009note}).

For the two-element Boolean algebra $\mathcal{B}=\{0,1\}$, $\mathcal{B}_n$ denotes the set of all $n \times n$ matrices over $\mathcal{B}$.
Under the commutative Boolean operations ($1 + 1 = 1$, $0 + 0 = 0$, $1 + 0 = 1$, $1 \times 1 = 1$, $0 \times 0 = 0$, $1 \times 0 = 0$), matrix addition and multiplication are still well-defined in $\mathcal{B}_n$.
Throughout this paper, a matrix is Boolean unless otherwise mentioned.

We note that the adjacency matrix of $C^m(D)$ for a digraph $D$ is the matrix $A^*_{m}$ obtained from $A^m(A^T)^m$ by replacing each of the diagonal elements with $0$ where $A$ is the adjacency matrix of $D$.
To see why, we take two distinct vertices $u$ and $v$ of $D$
corresponding to the $i$th row and the $j$th row of $A$, respectively.
Then
\begin{tabbing}
\ \ \ \ \ \ \= $u$ and $v$ are adjacent in $C^m(D)$ \\
$\Leftrightarrow$ \>  $u$ and $v$ have an $m$-step common prey in $D$ \\
$\Leftrightarrow$ \> the inner product of the $i$th row and the $j$th row of $A^m$ is $1$\\
$\Leftrightarrow$ \>  the $(i,j)$-entry of $A^*_m$ is $1$.
\end{tabbing}
It is easy to check that $A^*_i=A^*_j$ if and only if $A^i(A^T)^i=A^j(A^T)^j$ for a $(0,1)$ Boolean matrix $A$ of order $n$ and any positive integers $i$ and $j$.
Therefore, for a digraph $D$ and its adjacency matrix $A$,
\begin{equation}\label{eq:iff}
C^i(D)=C^j(D) \Leftrightarrow A^i(A^T)^i=A^j(A^T)^j
\end{equation}
for any positive integers $i$ and $j$.

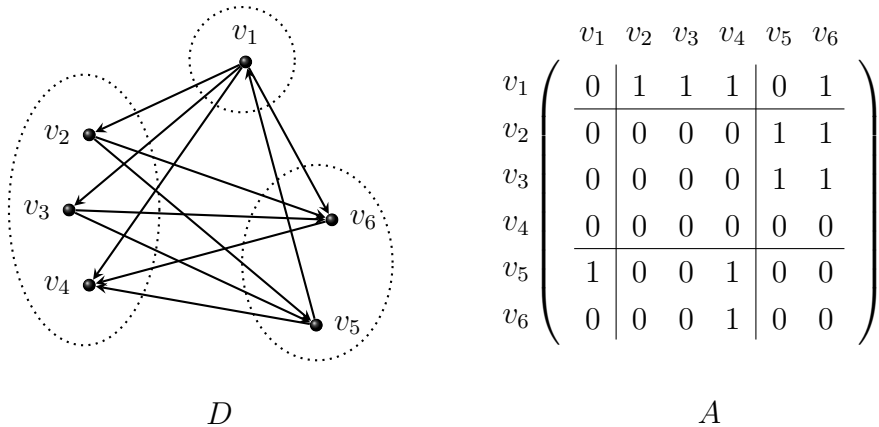
\begin{figure}
\begin{center}
\begin{minipage}[c]{.45\textwidth}
\begin{center}
\begin{tikzpicture}[auto,thick,scale=1]
    \tikzstyle{player}=[minimum size=5pt,inner sep=0pt,outer sep=0pt,ball color=black,circle]
    \tikzstyle{source}=[minimum size=5pt,inner sep=0pt,outer sep=0pt,ball color=black, circle]
    \tikzstyle{arc}=[minimum size=5pt,inner sep=1pt,outer sep=1pt, font=\footnotesize]
    \path (150:2cm)  node [player, label=left:$v_2$]  (2) {};
    \path (180:2cm)  node [player, label=left:$v_3$]  (3) {};
    \path (210:2cm)  node [player, label=left:$v_4$]  (4) {};
    \path (310:2cm)  node [player, label=right:$v_5$]  (5) {};
    \path (355:1.5cm)  node [player, label=right:$v_6$]  (6) {};
    \path (80:2cm)  node [player, label=above:$v_1$]  (1) {};
\draw[black,thick,-stealth] (1) - +(2);
\draw[black,thick,-stealth] (1) - +(3);
\draw[black,thick,-stealth] (1) - +(6);
\draw[black,thick,-stealth] (2) - +(5);
\draw[black,thick,-stealth] (2) - +(6);
\draw[black,thick,-stealth] (3) - +(5);
\draw[black,thick,-stealth] (3) - +(6);
\draw[black,thick,-stealth] (5) - +(1);
\draw[black,thick,-stealth] (6) - +(4);
\draw[black,thick,-stealth] (1) - +(4);
\draw[black,thick,-stealth] (5) - +(4);
\draw (0,-2.7) node{$D$};
\draw[dotted, thick] (-1.8, 0) ellipse (1 and 1.8);
\draw[dotted, thick] (1.3, -0.7) ellipse (1 and 1.3);
\draw[dotted, thick] (0.3, 2) ellipse (0.7 and 0.7);

    \end{tikzpicture}
\end{center}
\end{minipage}
\begin{minipage}[c]{.45\textwidth}
\begin{tikzpicture}
		[every matrix/.style={column sep = {1.5em,between origins},row sep={1.5em,between origins}}]
		\matrix (m) [matrix of math nodes,ampersand replacement=\&, left delimiter=(,right delimiter=)]
		{
			{0} \& {1} \& {1} \& {1} \& {0} \& {1} \\
			{0} \& {0} \& {0} \& {0} \& {1} \& {1} \\
			{0} \& {0} \& {0} \& {0} \& {1} \& {1} \\
			{0} \& {0} \& {0} \& {0} \& {0} \& {0} \\
			{1} \& {0} \& {0} \& {1} \& {0} \& {0} \\
			{0} \& {0} \& {0} \& {1} \& {0} \& {0} \\
		};
		
		\matrix (m2) [matrix of nodes,ampersand replacement=\&, left=0.4em of m]
		{
			{$v_1$} \\
			{$v_2$} \\
			{$v_3$} \\
			{$v_4$} \\
			{$v_5$} \\
			{$v_6$} \\
		};
		
		\matrix (m3) [matrix of nodes,ampersand replacement=\&, above=-0.4em of m]
		{
			{$v_1$} \& {$v_2$} \& {$v_3$} \& {$v_4$} \& {$v_5$} \& {$v_6$} \\
		};
\draw (0,-2.8) node{$A$};
\draw[-] ($0.5*(m-1-1.north east)+0.5*(m-1-2.north west)$) -- ($0.5*(m-6-1.south east)+0.5*(m-6-2.south west)$);
\draw[-] ($0.5*(m-1-1.south west)+0.5*(m-2-1.north west)$) -- ($0.5*(m-1-6.south east)+0.5*(m-2-6.north east)$);
\draw[-] ($0.5*(m-1-4.north east)+0.5*(m-1-5.north west)$) -- ($0.5*(m-6-4.south east)+0.5*(m-6-5.south west)$);
\draw[-] ($0.5*(m-4-1.south west)+0.5*(m-5-1.north west)$) -- ($0.5*(m-4-6.south east)+0.5*(m-5-6.north east)$);
\end{tikzpicture}
\end{minipage}
\end{center}
\caption{A tripartite tournament $D$ and its adjacency matrix $A$ where $V_1 =\{v_1\}$, $V_2 =\{v_2,v_3,v_4\}$, and $V_3=\{v_5, v_6\}$ are the partite sets of $D$.}
\label{fig:multipartitematrix}
\end{figure}

The greatest common divisor of the lengths of the closed directed walks of a strongly connected nontrivial digraph $D$ is called the {\it index of imprimitivity} of $D$ and denoted by $\kappa(D)$.
If $\kappa(D)=1$, then $D$ is said to be {\em primitive}.
The vertex set of $D$ can be partitioned into $\kappa(D)$ nonempty subsets $U_1, U_2, \ldots, U_{\kappa(D)}$ where each arc of $D$ goes out from $U_i$ and enters $U_{i+1}$ (identify $U_{\kappa(D)+1}$ with $U_1$) for some $i \in \{1 , \ldots, \kappa(D)\}$ (see \cite{brualdi1991combinatorial}).
We call the sets $U_1, U_2, \ldots, U_{\kappa(D)}$ the {\it sets of imprimitivity} of $D$.

A \emph{$k$-partite tournament} is an orientation of a complete $k$-partite graph for a positive integer $k$.
We call a $k$-partite tournament a {\em multipartite tournament} if $k \ge 2$.

Jung~{\em et al.}~\cite{jung2023limit} showed the following.
\statement{1234}
{
The index of imprimitivity of a strongly connected multipartite tournament can only be $1$, $2$, $3$, or $4$.
}
They also presented the sets of imprimitivity of a strongly connected multipartite tournament $D$ according to $\kappa(D)$ as follows.

\begin{Prop}[\cite{jung2023limit}]\label{cor:kappa}
Let $D$ be a strongly connected $k$-partite tournament with $k$-partition $(V_1,V_2,\ldots, V_k)$ for an integer $k \ge 2$.
Then the following are true:
\begin{itemize}
  \item[(i)] if $\kappa(D)=1$, then $k \ge 3$;
  \item[(ii)] if $\kappa(D)=2$, then $k=2$ and the sets of imprimitivity of $D$ are $V_1$ and $V_2$;
  \item[(iii)] if $\kappa(D)=3$, then $k=3$ and the sets of imprimitivity of $D$ are $V_1$, $V_2$, and $V_3$;
  \item[(iv)] if $\kappa(D)=4$, then $k=2$ and there exists a partition $\{X_i, Y_i\}$ of $V_i$ for each $i=1,2$ such that $X_1$, $Y_1$, $X_2$, and $Y_2$ are the sets of imprimitivity of $D$.
\end{itemize}
\end{Prop}

Suppose that $D$ is a weakly connected digraph.
Then the strong components of $D$ can be arranged as $Q_1,\ldots, Q_s$ so that there are no arcs going from $Q_i$ to $Q_j$ whenever $i > j$.
We call the strong components arranged in such a way
\textit{ordered strong components} of $D$.
A nontrivial component $Q_t$ is termed as the {\it last nontrivial component} of $D$ if either (i) $t=s$ or (ii) $t < s$ and $Q_i$ is trivial for each $t< i \le s$.
We denote the sets of imprimitivity of $Q_i$ by
\[
U_1^{(i)}, U_2^{(i)}, \ldots, U_{\kappa(Q_i)}^{(i)}
\]
for each integer $i=1,\ldots,s$.
Throughout this paper, we let $D_{\alpha \sim \beta}$ denote the subdigraph induced by $\bigcup_{i=\alpha}^{\beta}{V(Q_i)}$ if $\alpha \le \beta$ and be ignored if $\alpha > \beta$.

The adjacency matrix of a $k$-partite tournament for $k \ge 2$ can be represented as a block matrix $A$ with blocks $A_{ij}$ for $1 \le i,j \le k$ such that $A_{ii}$ is a zero matrix and $A_{ij}+A_{ji}^T$ is a matrix with all elements $1$ but not both corresponding elements of $A_{ij}$ and $A_{ji}^T$ equal to $1$ for $i \neq j$ (see $D$ and $A$ in Figure~\ref{fig:multipartitematrix} for an example).
In this paper, we determine the limit of the matrix sequence $\{A^m(A^T)^m\}_{m=1}^{\infty}$ in the case of convergence.
Further, in the event of divergence, we specify how the matrix sequence $\{A^m(A^T)^m\}_{m=1}^{\infty}$ changes periodically depending on the value of $m$.
Our results extend the work of Jung~{\em et al.}~\cite{jung2023limit} which addresses the case of the last strong component of its digraph being nontrivial, thereby completing the convergence analysis of $\{A^m(A^T)^m\}_{m=1}^{\infty}$.
To this end, we investigate the graph sequence $\{C^m(D)\}_{m=1}^{\infty}$ for the digraph $D$ of $A$, which is a multipartite tournament, based upon the observation \eqref{eq:iff}.

Given a $k$-partite tournament $D$ with ordered strong components $Q_1,\ldots, Q_s$ for some positive integers $k$ and $s$, let $Q_t$ be the last nontrivial component.
Then $Q_t$ is a strongly connected $l$-partite tournament for some positive integer $l$.
Now, by Proposition~\ref{cor:kappa}, if $\kappa(Q_t)=2$ (resp.\ $\kappa(Q_t)=3$), then $l=2$ (resp.\ $l=3$) and the sets of imprimitivity of $Q_t$ are the partite sets of $Q_t$; if $\kappa(Q_t)=4$, then $l=2$ and
each partite set of $Q_t$ can be partitioned into two subsets, which results in four subsets that make up the sets of imprimitivity of $Q_t$.
We also note that each partite set of $Q_t$ is the intersection of $Q_t$ and a partite set of $D$.
Based on these observations, we may assume, unless otherwise mentioned, that the partite sets of $D$ are labeled as $(V_1,\ldots,V_k)$ so that the sets of imprimitivity of $Q_t$ are represented in the following way:
\begin{equation}\label{eqn:bi}
  V_i \cap V(Q_t)=
  \begin{cases}
    U_{i}^{(t)}, & \mbox{if } \kappa(Q_t)=2 \mbox{ or } 3;\\
    U_{i}^{(t)} \cup U_{i+2}^{(t)}, &  \mbox{if } \kappa(Q_t)=4,
  \end{cases}
\end{equation}
for each $i=1,2,3$ if $\kappa(Q_t)=3$ and for each $i=1,2$ if $\kappa(Q_t)=2$ or $4$.

For vertex sets $X$ and $Y$ of a multipartite tournament $D$, we say that $X$ is {\it partite-related to} $Y$ if $X \subseteq W$ and $Y \subseteq W$ for some partite set $W$ of $D$.

Given a multipartite tournament $D$ with ordered strong components $Q_1,\ldots, Q_s$, let $X$ and $Y$ be partite sets  of $Q_i$ and $Q_j$, respectively, for some $1 \le i < j \le s$ such that $X$ is not partite-related to $Y$.
Then $X \to Y$.

Theorem~\ref{Thm:diverse} completely characterizes the convergence of $\{C^m(D)\}_{m=1}^{\infty}$ for a multipartite tournament $D$ based on the last nontrivial component of $D$.

\begin{Thm}\label{Thm:diverse}
  Let $D$ be a $k$-partite tournament with $k$-partition $(V_1,\ldots,V_k)$ and ordered strong components $Q_1,\ldots, Q_s$.
  Then the graph sequence $\{C^m(D)\}_{m=1}^{\infty}$ diverges if and only if (a) $s \ge 2 $, (b) $D$ has a nontrivial strong component, (c) for the last nontrivial component $Q_t$, $t < s$, and (d) one of the following is true:
  \begin{itemize}
    \item[(i)] $\kappa(Q_t)=3$ and $\bigcup_{i=t+1}^{s}{V(Q_i)} \subseteq V_j$ for some $j \in \{1, 2, 3\}$;
    \item[(ii)] $\kappa(Q_t)=3$ and there is an integer $\alpha \in \{t+1, \ldots, s-1\}$ such that $V(D_{t+1 \sim \alpha })$ and $V(D_{\alpha+1 \sim s})$ are partite-related to $U_{j}^{(t)}$ and $U_{j+1}^{(t)}$, respectively, for some integer $j \in \{1,2,3\}$ (identifying $U_{4}^{(t)}$ with $U_{1}^{(t)}$);
    \item[(iii)] $\kappa(Q_t)=4$ and $\bigcup_{i=t+1}^{s}{V(Q_i)} \subseteq V_j$ for some $j \in \{1, 2\}$.
    \end{itemize}
\end{Thm}

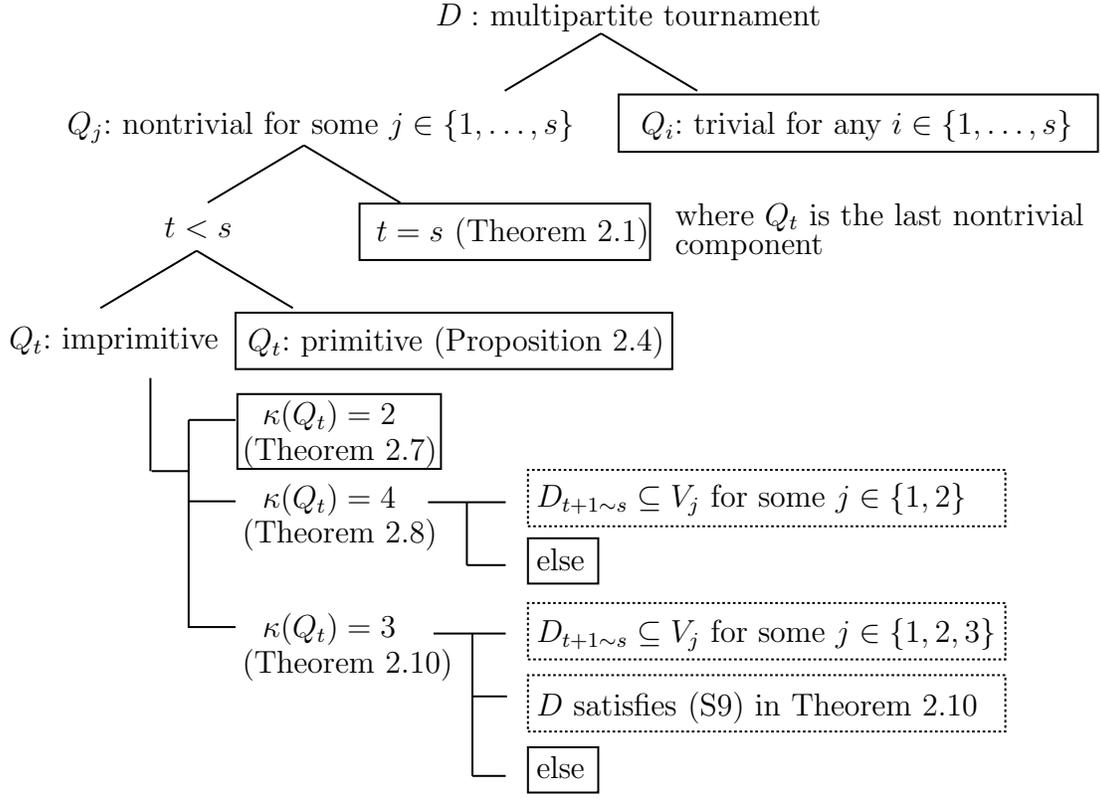
\begin{figure}[h]
  \centering
  \tikzset{every picture/.style={line width=0.75pt}} 

\begin{tikzpicture}[x=0.7pt,y=0.7pt,yscale=-1,xscale=1]

\draw    (269.5,51) -- (321.5,20) ;
\draw    (321.5,20) -- (373.5,51) ;

\draw    (109,111.5) -- (161,80.5) ;
\draw    (161,80.5) -- (213,111.5) ;

\draw    (51,168.5) -- (103,137.5) ;
\draw    (103,137.5) -- (155,168.5) ;

\draw    (98.67,229) -- (98.67,341.5) ;
\draw    (98.67,229) -- (124,229) ;
\draw    (78,206.33) -- (78,256.67) ;
\draw    (78.67,256.67) -- (98.67,256.67) ;
\draw    (98.67,273) -- (124,273) ;
\draw    (98.67,341) -- (124,341) ;
\draw    (228,273.5) -- (270,273.5) ;
\draw    (249,273.5) -- (249,307.5) ;
\draw    (249,307.5) -- (270,307.5) ;
\draw    (231,344.5) -- (270,344.5) ;
\draw    (252,344.5) -- (252,421.5) ;
\draw    (252,378.5) -- (270.67,378.5) ;
\draw    (252,421.5) -- (270,421.5) ;
\draw  [dash pattern={on 1pt off 1pt}] (282,256) -- (540,256) -- (540,286.5) -- (282,286.5) -- cycle ;
\draw  [dash pattern={on 1pt off 1pt}] (282,328) -- (540,328) -- (540,358.5) -- (282,358.5) -- cycle ;
\draw  [dash pattern={on 1pt off 1pt}] (282,367) -- (540,367) -- (540,397.5) -- (282,397.5) -- cycle ;
\draw   (191,112) -- (348,112) -- (348,142.5) -- (191,142.5) -- cycle ;
\draw   (331,53) -- (590,53) -- (590,84) -- (331,84) -- cycle ;
\draw   (124,171) -- (360,171) -- (360,201.5) -- (124,201.5) -- cycle ;
\draw   (125,215) -- (235,215) -- (235,255.5) -- (125,255.5) -- cycle ;
\draw   (282,293) -- (320,293) -- (320,317.5) -- (282,317.5) -- cycle ;
\draw   (282,406) -- (320,406) -- (320,430.5) -- (282,430.5) -- cycle ;

\draw (198.5,118) node [anchor=north west][inner sep=0.75pt]    {$t=s$ (Theorem~\ref{thm:last})};
\draw (360,110) node [anchor=north west][inner sep=0.75pt]    {where $Q_t$ is the last nontrivial};
\draw (360,127) node [anchor=north west][inner sep=0.75pt]    {component};
\draw (84,118) node [anchor=north west][inner sep=0.75pt]    {$t<s$};
\draw (129,177.07) node [anchor=north west][inner sep=0.75pt]    {$Q_{t}$: primitive (Proposition~\ref{prop:qt_primitive})};
\draw (0,177.07) node [anchor=north west][inner sep=0.75pt]    {$Q_{t}$: imprimitive};
\draw (137,217) node [anchor=north west][inner sep=0.75pt]    {$\kappa ( Q_{t}) =2$};
\draw (137,262) node [anchor=north west][inner sep=0.75pt]    {$\kappa ( Q_{t}) =4$};
\draw (137,332) node [anchor=north west][inner sep=0.75pt]    {$\kappa ( Q_{t}) =3$};

\draw (126,281) node [anchor=north west][inner sep=0.75pt]   [align=left] {(Theorem~\ref{thm:4})};
\draw (126,351) node [anchor=north west][inner sep=0.75pt]   [align=left] {(Theorem~\ref{thm:3})};
\draw (126,236) node [anchor=north west][inner sep=0.75pt]   [align=left] {(Theorem~\ref{thm:k=2})};

\draw (31.5,60) node [anchor=north west][inner sep=0.75pt]    {$Q_{j}$: nontrivial for some $j \in \{1,\dotsc ,s\}$};
\draw (285,261.4) node [anchor=north west][inner sep=0.75pt]    {$D_{t+1 \sim  s} \subseteq V_{j}$ for some $j \in \{1,2\}$};
\draw (285,334.4) node [anchor=north west][inner sep=0.75pt]    {$D_{t+1 \sim s} \subseteq V_{j}$ for some $j \in  \{1,2,3\}$};
\draw (285,297) node [anchor=north west][inner sep=0.75pt]   [align=left] {else};
\draw (285,409) node [anchor=north west][inner sep=0.75pt]   [align=left] {else};
\draw (230.5,1.9) node [anchor=north west][inner sep=0.75pt]    {$D:$ multipartite tournament};
\draw (285,374.4) node [anchor=north west][inner sep=0.75pt]    {$D$ satisfies \eqref{dag} in Theorem~\ref{thm:3}};
\draw (341.5,60) node [anchor=north west][inner sep=0.75pt]    {$Q_{i}$: trivial for any $i\in \{1,\dotsc ,s\}$};
\end{tikzpicture}
  \caption{Solid-line boxes represent the converging cases, while dashed-line boxes represent the diverging cases.}\label{fig:diverse}
\end{figure}

Theorem~\ref{Thm:diverse} will be proven through the entire paper.
The overview of this paper is given in Figure~\ref{fig:diverse}.

Furthermore, not only do we determine the limit in the case of convergence, but also in the event of divergence, we specify how $C^m(D)$ changes periodically depending on the value of $m$.
Refer to Theorems~\ref{thm:last}, \ref{thm:k=2}, \ref{thm:4}, \ref{thm:3} and Proposition~\ref{prop:qt_primitive}.
Our results extend the work of Jung~{\em et al.}~\cite{jung2023limit} which addresses the case of the last strong component being nontrivial  (see Theorem~\ref{thm:last}), thereby completing the convergence analysis of  $\{C^m(D)\}$ for a multipartite tournament $D$.

By \eqref{eq:iff}, the convergence of $\{C^m(D)\}_{m=1}^{\infty}$ is consistent with that of the matrix sequence $\{A^m(A^T)^m\}_{m=1}^{\infty}$ for a digraph $D$ and its adjacency matrix $A$.

\begin{figure}[h]
\begin{center}
\begin{minipage}[c]{.45\textwidth}
\begin{center}
\begin{tikzpicture}[auto,thick,scale=1]
    \tikzstyle{player}=[minimum size=5pt,inner sep=0pt,outer sep=0pt,ball color=black,circle]
    \tikzstyle{source}=[minimum size=5pt,inner sep=0pt,outer sep=0pt,ball color=black, circle]
    \tikzstyle{arc}=[minimum size=5pt,inner sep=1pt,outer sep=1pt, font=\footnotesize]
    \path (110:2cm)  node [player, label=above:$v_3$]  (6) {};
    \path (180:2cm)  node [player, label=left:$v_2$]  (1) {};
    \path (210:2cm)  node [player, label=left:$v_1$]  (2) {};
    \path (330:2cm)  node [player, label=right:$v_6$]  (3) {};
    \path (0:2cm)  node [player, label=right:$v_5$]  (4) {};
    \path (70:2cm)  node [player, label=above:$v_4$]  (5) {};
\draw[black,thick,-stealth] (1) - +(3);
\draw[black,thick,-stealth] (1) - +(4);
\draw[black,thick,-stealth] (1) - +(6);
\draw[black,thick,-stealth] (2) - +(3);
\draw[black,thick,-stealth] (2) - +(4);
\draw[black,thick,-stealth] (2) - +(6);
\draw[black,thick,-stealth] (3) - +(5);
\draw[black,thick,-stealth] (3) - +(6);
\draw[black,thick,-stealth] (4) - +(6);
\draw[black,thick,-stealth] (5) - +(1);
\draw[black,thick,-stealth] (5) - +(2);
\draw[black,thick,-stealth] (5) - +(4);

\draw (-3.2,-0.4) node{$V_1$};
\draw (3.2,-0.4) node{$V_2$};
\draw (0,3.3) node{$V_3$};
\draw[dotted, thick] (-1.8, -0.4) ellipse (1 and 1.3);
\draw[dotted, thick] (1.8, -0.4) ellipse (1 and 1.3);
\draw[dotted, thick] (0, 2) ellipse (1.3 and 1);

    \end{tikzpicture}
    $D_1$
\end{center}
\end{minipage}
\phantom{dd}
\begin{minipage}[c]{.45\textwidth}

\begin{tikzpicture}[auto,thick,scale=1]
    \tikzstyle{player}=[minimum size=5pt,inner sep=0pt,outer sep=0pt,ball color=black,circle]
    \tikzstyle{source}=[minimum size=5pt,inner sep=0pt,outer sep=0pt,ball color=black, circle]
    \tikzstyle{arc}=[minimum size=5pt,inner sep=1pt,outer sep=1pt, font=\footnotesize]
    \path (180:4cm)  node [player, label=above:$v_1$]  (1) {};
    \path (160:3.5cm)  node [player, label=left:$v_2$]  (2) {};
    \path (160:2cm)  node [player, label=above:$v_6$]  (6) {};
    \path (0:0cm)  node [player, label=below:$v_5$]  (5) {};
    \path (200:3cm)  node [player, label=below:$v_4$]  (4) {};
    \path (0:2cm)  node [player, label=below:$v_3$]  (3) {};
\draw[black,thick,-stealth] (1) - +(5);
\draw[black,thick,-stealth] (1) - +(6);
\draw[black,thick,-stealth] (2) - +(6);
\draw[black,thick,-stealth] (4) - +(1);
\draw[black,thick,-stealth] (4) - +(2);
\draw[black,thick,-stealth] (4) - +(5);
\draw[black,thick,-stealth] (5) - +(3);
\draw[black,thick,-stealth] (6) - +(3);
\draw[black,thick,-stealth] (6) - +(4);
 \path (2) edge [black, bend left=20,thick,-stealth] (3);
\path (1) edge [black, bend right=20,thick,-stealth] (3);
\path (2) edge [black, bend right=15,thick,-stealth] (5);

\draw[rounded corners, dash pattern={on 1pt off 1pt}] (-4.5, -1.7) rectangle (-1.3, 1.7) {};
\draw[rounded corners, dash pattern={on 1pt off 1pt}] (-0.5, -0.7) rectangle (0.7, 0.5) {};
\draw[rounded corners, dash pattern={on 1pt off 1pt}] (1.5, -0.7) rectangle (2.7, 0.5) {};
\draw (-2.8,-2) node{$Q_1$};
\draw (0.1,-1) node{$Q_2$};
\draw (2.2,-1) node{$Q_3$};
    \end{tikzpicture}

\end{minipage}
\end{center}
\begin{center}
\begin{minipage}[c]{.45\textwidth}
\begin{center}
\begin{tikzpicture}[auto,thick,scale=1]
    \tikzstyle{player}=[minimum size=5pt,inner sep=0pt,outer sep=0pt,ball color=black,circle]
    \tikzstyle{source}=[minimum size=5pt,inner sep=0pt,outer sep=0pt,ball color=black, circle]
    \tikzstyle{arc}=[minimum size=5pt,inner sep=1pt,outer sep=1pt, font=\footnotesize]
    \path (150:2cm)  node [player, label=above:$v_3$]  (6) {};
    \path (180:2cm)  node [player, label=left:$v_2$]  (1) {};
    \path (210:2cm)  node [player, label=left:$v_1$]  (2) {};
    \path (330:2cm)  node [player, label=right:$v_6$]  (3) {};
    \path (0:2cm)  node [player, label=right:$v_5$]  (4) {};
    \path (75:2.3cm)  node [player, label=above:$v_4$]  (5) {};
\draw[black,thick,-stealth] (1) - +(3);
\draw[black,thick,-stealth] (1) - +(4);
\draw[black,thick,-stealth] (2) - +(3);
\draw[black,thick,-stealth] (2) - +(4);
\draw[black,thick,-stealth] (3) - +(5);
\draw[black,thick,-stealth] (3) - +(6);
\draw[black,thick,-stealth] (4) - +(6);
\draw[black,thick,-stealth] (5) - +(1);
\draw[black,thick,-stealth] (5) - +(2);
\draw[black,thick,-stealth] (5) - +(4);
\draw[black,thick,-stealth] (5) - +(6);

\draw[dotted, thick] (-1.8, 0) ellipse (1.2 and 2);
\draw[dotted, thick] (1.8, -0.4) ellipse (1 and 1.3);
\draw[dotted, thick] (0.5, 2.2) ellipse (1 and 1);

\draw (-3.4,0.2) node{$V_1$};
\draw (3.2,-0.4) node{$V_2$};
\draw (0.8,3.5) node{$V_3$};

    \end{tikzpicture}
$D_2$
\end{center}
\end{minipage}
\phantom{dd}
\begin{minipage}[c]{.45\textwidth}
\begin{tikzpicture}[auto,thick,scale=1]
    \tikzstyle{player}=[minimum size=5pt,inner sep=0pt,outer sep=0pt,ball color=black,circle]
    \tikzstyle{source}=[minimum size=5pt,inner sep=0pt,outer sep=0pt,ball color=black, circle]
    \tikzstyle{arc}=[minimum size=5pt,inner sep=1pt,outer sep=1pt, font=\footnotesize]
    \path (180:4cm)  node [player, label=above:$v_1$]  (1) {};
    \path (160:3.5cm)  node [player, label=left:$v_2$]  (2) {};
    \path (160:2cm)  node [player, label=above:$v_6$]  (6) {};
    \path (0:0cm)  node [player, label=below:$v_5$]  (5) {};
    \path (200:3cm)  node [player, label=below:$v_4$]  (4) {};
    \path (0:2cm)  node [player, label=below:$v_3$]  (3) {};
\draw[black,thick,-stealth] (1) - +(5);
\draw[black,thick,-stealth] (1) - +(6);
\draw[black,thick,-stealth] (2) - +(6);
\draw[black,thick,-stealth] (4) - +(1);
\draw[black,thick,-stealth] (4) - +(2);
\draw[black,thick,-stealth] (4) - +(5);
\draw[black,thick,-stealth] (5) - +(3);
\draw[black,thick,-stealth] (6) - +(3);
\draw[black,thick,-stealth] (6) - +(4);
\path (2) edge [black, bend right=15,thick,-stealth] (5);
\path (4) edge [black, bend right=10,thick,-stealth] (3);

\draw[rounded corners, dash pattern={on 1pt off 1pt}] (-4.5, -1.7) rectangle (-1.3, 1.7) {};
\draw[rounded corners, dash pattern={on 1pt off 1pt}] (-0.5, -0.7) rectangle (0.7, 0.5) {};
\draw[rounded corners, dash pattern={on 1pt off 1pt}] (1.5, -0.7) rectangle (2.7, 0.5) {};
\draw (-2.8,-2) node{$Q_1$};
\draw (0.1,-1) node{$Q_2$};
\draw (2.2,-1) node{$Q_3$};
    \end{tikzpicture}
\end{minipage}

\end{center}
\caption{Tripartite tournaments $D_1$ and $D_2$ and their ordered strong components $Q_1, Q_2, Q_3$.}
\label{fig:multipartitematrix2}
\end{figure}
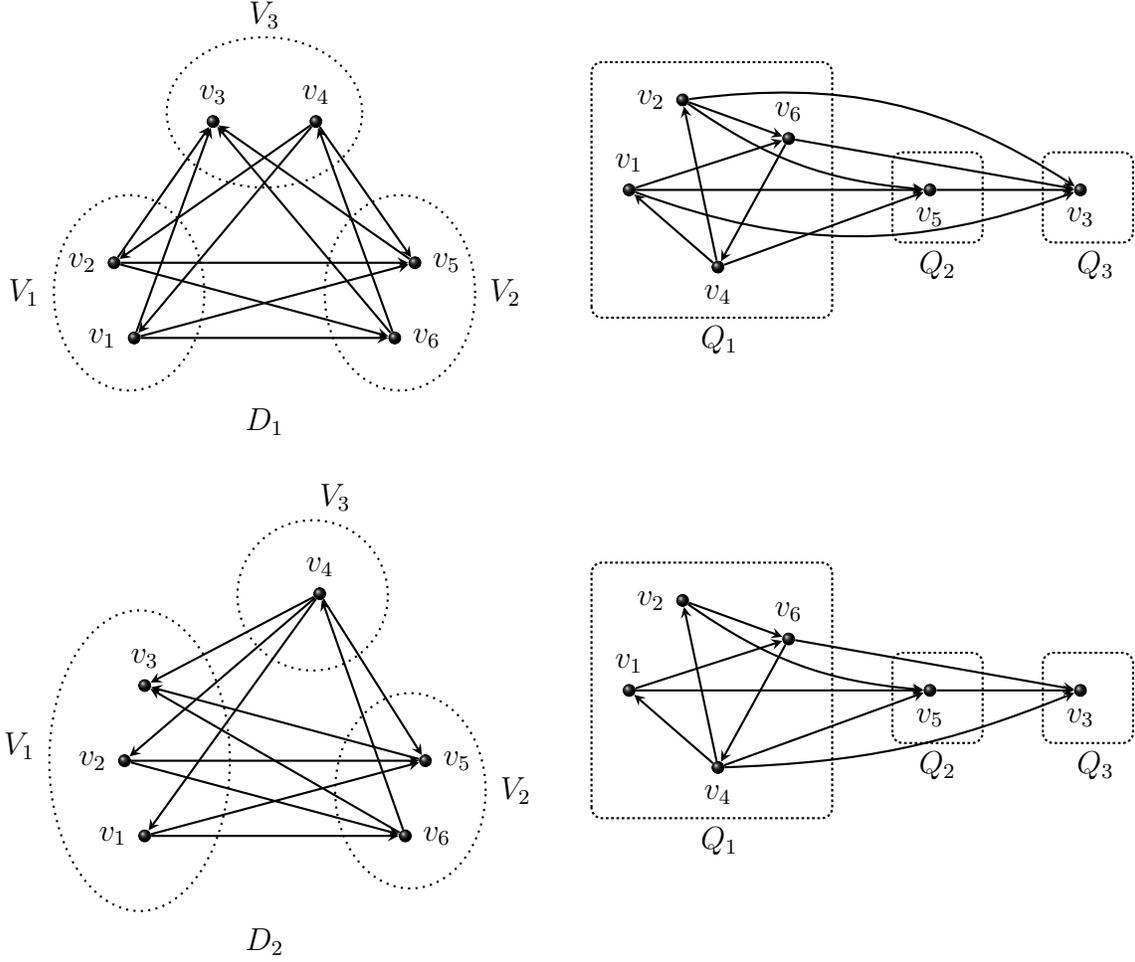

\begin{Ex}
Consider the matrix sequences $\{A_1^m (A_1^T)^m\}_{m=1}^{\infty}$ and $\{A_2^m (A_2^T)^m\}_{m=1}^{\infty}$ where
  \begin{equation*}
    A_1 =
\left(    \begin{matrix}
    0 & 0 & 1 & 0 & 1 & 1 \\
    0 & 0 & 1 & 0 & 1 & 1 \\
    0 & 0 & 0 & 0 & 0 & 0 \\
    1 & 1 & 0 & 0 & 1 & 0 \\
    0 & 0 & 1 & 0 & 0 & 0 \\
    0 & 0 & 1 & 1 & 0 & 0
    \end{matrix} \right)
, \quad    A_2 =
\left(    \begin{matrix}
    0 & 0 & 0 & 0 & 1 & 1 \\
    0 & 0 & 0 & 0 & 1 & 1 \\
    0 & 0 & 0 & 0 & 0 & 0 \\
    1 & 1 & 1 & 0 & 1 & 0 \\
    0 & 0 & 1 & 0 & 0 & 0 \\
    0 & 0 & 1 & 1 & 0 & 0
    \end{matrix} \right).
  \end{equation*}
To determine whether or not they converge, we consider the digraphs $D_1$ and $D_2$ of $A_1$ and $A_2$, respectively.
   See Figure~\ref{fig:multipartitematrix2} for an illustration  where $v_i$ is the vertex corresponding to the $i$th row for $i =1, 2, \ldots, 6$.
For the ordered strong components $Q_1, Q_2, Q_3$, $Q_1$ is the last nontrivial component in each of digraphs $D_1$ and $D_2$.
In addition, $\kappa(Q_1)=3$, $U_{1}^{(1)} = \{v_1, v_2\}$, $U_{2}^{(1)}=\{v_6\}$, and $U_{3}^{(1)}=\{v_4\}$ for $D_1$ and $D_2$.

In the case of $D_1$, $V(Q_2)$ and $V(Q_3)$ are partite-related to $U_{2}^{(1)}$ and $U_{3}^{(1)}$, respectively.
Thus \eqref{dag} in Theorem~\ref{thm:3} is satisfied and so the graph sequence $\{C^m(D_1)\}_{m=1}^{\infty}$ diverges by Theorem~\ref{Thm:diverse}.
Therefore the matrix sequences $\{A_1^m (A_1^T)^m\}_{m=1}^{\infty}$ diverges by \eqref{eq:iff}.
Further, $A_1^{3m+2} (A_1^T)^{3m+2}$, $A_1^{3m+3} (A_1^T)^{3m+3}$, and $A_1^{3m+4} (A_1^T)^{3m+4}$ equal the first, the second, and the third matrix below, respectively, for any nonnegative integer $m$.

\begin{equation*}
  \left(    \begin{matrix}
    1 & 1 & 0 & 1 & 0 & 0 \\
    1 & 1 & 0 & 1 & 0 & 0 \\
    0 & 0 & 0 & 0 & 0 & 0 \\
    1 & 1 & 0 & 1 & 0 & 1 \\
    0 & 0 & 0 & 0 & 0 & 0 \\
    0 & 0 & 0 & 1 & 0 & 1
    \end{matrix} \right),
  \quad   \left(    \begin{matrix}
    1 & 1 & 0 & 0 & 0 & 1 \\
    1 & 1 & 0 & 0 & 0 & 1 \\
    0 & 0 & 0 & 0 & 0 & 0 \\
    0 & 0 & 0 & 1 & 0 & 1 \\
    0 & 0 & 0 & 0 & 0 & 0 \\
    1 & 1 & 0 & 1 & 0 & 1
    \end{matrix} \right),
  \quad   \left(    \begin{matrix}
    1 & 1 & 0 & 1 & 0 & 1 \\
    1 & 1 & 0 & 1 & 0 & 1 \\
    0 & 0 & 0 & 0 & 0 & 0 \\
    1 & 1 & 0 & 1 & 0 & 0 \\
    0 & 0 & 0 & 0 & 0 & 0 \\
    1 & 1 & 0 & 0 & 0 & 1
    \end{matrix} \right)
\end{equation*}

In the case of $D_2$, $V(Q_2)$ and $V(Q_3)$ are partite-related to $U_{2}^{(1)}$ and $U_{1}^{(1)}$, respectively.
Thus $D_{2 \sim 3} \not\subseteq V_j$ for any $j \in \{1,2,3\}$  and \eqref{dag} in Theorem~\ref{thm:3} is not satisfied.
Therefore the graph sequence $\{C^m(D_2)\}_{m=1}^{\infty}$ converges by Theorem~\ref{Thm:diverse} and so the matrix sequences $\{A_2^m (A_2^T)^m\}_{m=1}^{\infty}$ converges by \eqref{eq:iff}.
Further, its limit is equal to
\begin{equation*}
  \left(    \begin{matrix}
    1 & 1 & 0 & 1 & 0 & 1 \\
    1 & 1 & 0 & 1 & 0 & 1 \\
    0 & 0 & 0 & 0 & 0 & 0 \\
    1 & 1 & 0 & 1 & 0 & 1 \\
    0 & 0 & 0 & 0 & 0 & 0 \\
    1 & 1 & 0 & 1 & 0 & 1
    \end{matrix} \right).
\end{equation*}
\end{Ex}

\section{A Proof of Theorem~\ref{Thm:diverse}}

Let $G_1$ and $G_2$ be vertex-disjoint graphs.
We call the graph having the vertex set $V(G_1)\cup V(G_2)$ and the edge set $E(G_1)\cup E(G_2)$ the {\it union} of $G_1$ and $G_2$ and denote it by $G_1 \cup G_2$.
We call the graph having the vertex set $V(G_1)\cup V(G_2)$ and the edge set $E(G_1)\cup E(G_2) \cup \{uv \mid u \in V(G_1), v \in V(G_2)\}$ the {\it join} of $G_1$ and $G_2$ and denote it by $G_1 \vee G_2$ .

Given a vertex set $Z$ of a digraph $D$, we denote by $K[Z]$ (resp.\ $I[Z]$) the complete graph (resp.\ the empty graph) with vertex set $Z$.
In particular, we write $K[V(D)]$ and $I[V(D)]$ as $K[D]$ and $I[D]$, respectively.

In this paper, we study the convergence of the graph sequence $\{C^m(D)\}_{m=1}^{\infty}$ for a multipartite tournament $D$ in terms of ordered strong components of $D$.
Given a digraph $D$, if $D$ has no nontrivial strong components, then $\{C^m(D)\}_{m=1}^{\infty}$ converges to an empty graph.
Therefore we only consider the multipartite tournament with at least one nontrivial strong component.
The limit of $\{C^m(D)\}_{m=1}^{\infty}$ in the case where the last strong component is nontrivial is well computed by Jung~{\em et al.}~\cite{jung2023limit} and obtained as follows.

\begin{figure}
  \begin{center}
  \begin{tabular}{cc}
\begin{tikzpicture}[auto,thick,scale=1]
    \tikzstyle{player}=[minimum size=5pt,inner sep=1pt,outer sep=0pt,circle,draw=black]
    \tikzstyle{player1}=[minimum size=5pt,inner sep=1pt,outer sep=3pt,rectangle,draw=white]
    \tikzstyle{source}=[minimum size=5pt,inner sep=0pt,outer sep=0pt,ball color=black, circle]
    \tikzstyle{arc}=[minimum size=5pt,inner sep=1pt,outer sep=1pt, font=\footnotesize]
    \path (90:2cm)     node [player1]  (a) {\phantom{d}$K^{(1)}$\phantom{d}};
    \path (0:2cm)       node [player1] (b){\phantom{dd}$K^{(3)}$\phantom{dd}};
    \path (180:2cm)       node [player1] (c){\phantom{dd}$K^{(2)}$\phantom{dd}};
    \draw[black,thick,-] (a) -- ++(b);
    \draw[black,thick,-] (a) -- ++(c);

    \end{tikzpicture}
    &
\begin{tikzpicture}
\matrix [matrix of math nodes,left delimiter=(,right delimiter=),row sep=0.3cm,column sep=0.1cm] (m) {
      J^{(1)} & J &  J \\
      J & J^{(2)} &  O  \\
      J & O   &J^{(3)}\\ };
\end{tikzpicture}
\\
    $G_1$ & $M_1$\\
    &\\

    \begin{tikzpicture}[auto,thick,scale=0.8]
    \tikzstyle{player}=[minimum size=5pt,inner sep=1pt,outer sep=0pt,circle,draw=black]
    \tikzstyle{player1}=[minimum size=5pt,inner sep=1pt,outer sep=3pt,rectangle,draw=white]
    \tikzstyle{source}=[minimum size=5pt,inner sep=0pt,outer sep=0pt,ball color=black, circle]
    \tikzstyle{arc}=[minimum size=5pt,inner sep=1pt,outer sep=1pt, font=\footnotesize]
    \path (90:3cm)     node [player1]  (a) {$K^{(1)}$};
    \path (300:3.6cm)       node [player1] (b){$K^{(3)}$};
    \path (240:3.6cm)       node [player1] (c){$K^{(2)}$};
    \path (0:2cm)       node [player1] (d){$K^{(6)}$};
    \path (0:4cm)       node [player1] (e){$K^{(7)}$};
    \path (180:2cm)       node [player1] (f){$K^{(5)}$};
    \path (180:4cm)       node [player1] (g){$K^{(4)}$};

    \draw[black,thick,-] (a) -- ++(b);
    \draw[black,thick,-] (a) -- ++(c);
    \draw[black,thick,-] (b) -- ++(d);
    \draw[black,thick,-] (b) -- ++(e);
    \draw[black,thick,-] (c) -- ++(f);
    \draw[black,thick,-] (c) -- ++(g);
    \draw[black,thick,-] (a) -- ++(f);
    \draw[black,thick,-] (a) -- ++(g);
    \draw[black,thick,-] (a) -- ++(d);
    \draw[black,thick,-] (a) -- ++(e);

    \end{tikzpicture}
    & \begin{tikzpicture}
\matrix [matrix of math nodes,left delimiter=(,right delimiter=),row sep=0.1cm,column sep=-0.1cm] (m) {
      J^{(1)} & J & J  & J & J  & J & J  \\
      J & J^{(2)} & O  & J & J  & O & O  \\
      J & O   &J^{(3)} & O & O  & J & J  \\
      J & J & O & J^{(4)} & O & O & O \\
      J & J & O &O   &J^{(5)}& O & O  \\
      J &  O & J  &  O   & O  & J^{(6)} & O  \\
      J & O & J & O & O & O & J^{(7)}\\ };
\end{tikzpicture}\\
$G_2$ & $M_2$\\
    &\\

    \begin{tikzpicture}[auto,thick,scale=0.7]
    \tikzstyle{player}=[minimum size=5pt,inner sep=1pt,outer sep=0pt,circle,draw=black]
    \tikzstyle{player1}=[minimum size=5pt,inner sep=0pt,outer sep=3pt,rectangle,draw=white]
    \tikzstyle{source}=[minimum size=5pt,inner sep=0pt,outer sep=0pt,ball color=black, circle]
    \tikzstyle{arc}=[minimum size=5pt,inner sep=1pt,outer sep=1pt, font=\footnotesize]
       \path (40:5cm)  node [player1]  (a) {$K^{(4)}$};
    \path (90:1.5cm)     node [player1]  (b) {$K^{(3)}$};
    \path (140:5cm)  node [player1]  (c) {$K^{(2)}$};
    \path (340:3.6cm)   node [player1]  (d) {$K^{(7)}$};
    \path (270:2cm)       node [player1]  (e){$K^{(6)}$};
    \path (200:3.6cm)  node [player1]  (f){$K^{(5)}$};
    \path (90:5cm)  node [player1]  (g){$K^{(1)}$};
    \draw[black,thick,-] (a) -- ++(b);
    \draw[black,thick,-] (b) -- ++(c);
   \draw[black,thick,-] (a) -- ++(c);
    \draw[black,thick,-] (a) -- ++(e);
    \draw[black,thick,-] (c) -- ++(e);
    \draw[black,thick,-] (b) -- ++(f);
    \draw[black,thick,-] (b) -- ++(d);
    \draw[black,thick,-] (c) -- ++(f);
    \draw[black,thick,-] (a) -- ++(d);
    \draw[black,thick,-] (g) -- ++(a);
    \draw[black,thick,-] (g) -- ++(b);
    \draw[black,thick,-] (g) -- ++(c);
    \draw[black,thick,-] (g) -- ++(d);
   \path (g) edge [black, bend right=16,thick,-] (e);
    \draw[black,thick,-] (g) -- ++(f);
    \end{tikzpicture}
         &
        \begin{tikzpicture}
\matrix [matrix of math nodes,left delimiter=(,right delimiter=),row sep=0.1cm,column sep=-0.1cm] (m) {
      J^{(1)} & J &   J  & J & J   & J &J \\
      J & J^{(2)} &   J  & J & J   & J & O  \\
      J & J   &J^{(3)} & J & J   & O & J  \\
      J & J & J & J^{(4)} & O & J & J \\
      J & J    & J    &   O   &J^{(5)}& O & O  \\
      J &  J    & O    &  J   & O  & J^{(6)} &O  \\
      J & O & J & J & O & O & J^{(7)} \\ };
\end{tikzpicture}\\
        $G_3$ & $M_3$
\end{tabular}
\end{center}
\caption{The graphs $G_1$, $G_2$, and $G_3$ in Theorem~\ref{thm:last} (refer to \eqref{figure} for other notations).
}\label{fig:graphs}
\end{figure}
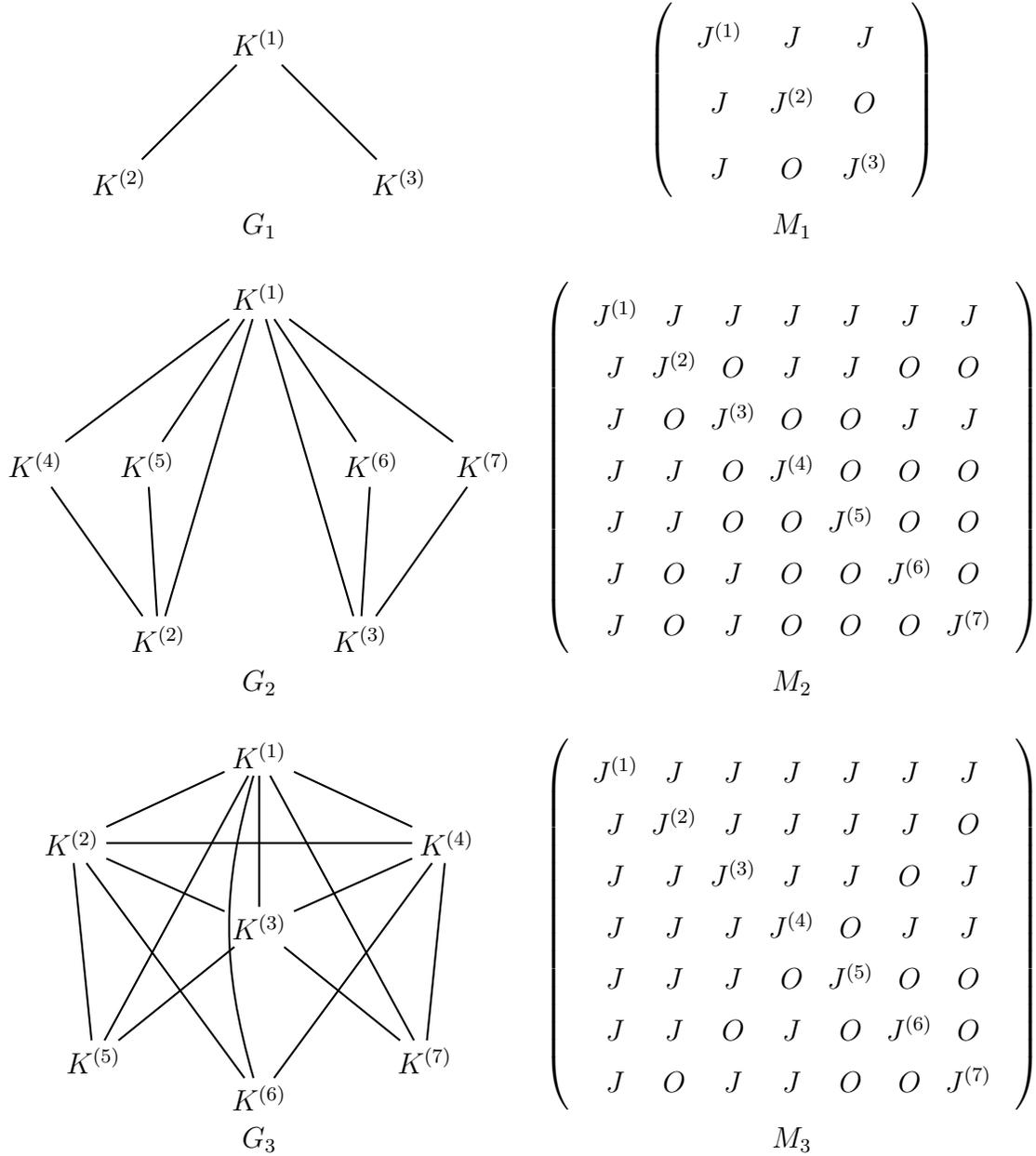

\begin{Thm}[\cite{jung2023limit}]\label{thm:last}
Let $D$ be a multipartite tournament with ordered strong components $Q_1,\ldots, Q_s$ for some integer $s \ge 1$.
If $Q_s$ is nontrivial, then the graph sequence $\{C^m(D)\}_{m=1}^{\infty}$ converges to a graph
\[
G \cong \begin{cases}
          K[D], & \mbox{if } \kappa(Q_s)=1; \\
          G_1, & \mbox{if } \kappa(Q_s)=2; \\
          G_3, & \mbox{if } \kappa(Q_s)=3; \\
          G_2, & \mbox{if } \kappa(Q_s)=4,
        \end{cases}
\]
 where $G_1$, $G_2$, and $G_3$ are the graphs given in Figure~\ref{fig:graphs}.
 \end{Thm}

In each graph given in Figure~\ref{fig:graphs} or  \ref{fig:graphs_for4} or \ref{fig:graphs_for3}, 
\statement{figure}{
\begin{itemize}
	\item $K^{(i)}$ either does not exist or stands for a clique;
	\item the line between two cliques $K^{(i)}$ and $K^{(j)}$ indicates that $V\left(K^{(i)}\right) \cup V\left(K^{(j)}\right)$ forms a clique while the absence of a line between $K^{(i)}$ and $K^{(j)}$ means that there are no edges joining a vertex in $K^{(i)}$ and a vertex in $K^{(j)}$;
	\item the matrix $M_i= A(G_i)+I$ where $A(G_i)$ is the adjacency matrix of $G_i$ for each $i=1,2,3,4,5$;
	\item $I$, $J$, and $O$ are an identity matrix, a matrix of all $1$s, and a zero matrix, respectively, of an appropriate order;
	\item the block matrix $J^{(i)}$ represents the clique $K^{(i)}$.
\end{itemize}
}

According to \eqref{1234}, Theorem~\ref{thm:last} completes computing the limit of the graph sequence $\{C^m(D)\}_{m=1}^{\infty}$ for a multipartite tournament $D$ having the nontrivial last strong component.
Consequently, it remains to consider the multipartite tournaments with at least one nontrivial strong component that is not the last strong component of $D$.
From now on, $Q_s$ is trivial for ordered strong components $Q_1,\ldots, Q_s$ of $D$ for some integer $s \ge 2$.

Let $D$ be a multipartite tournament with ordered strong components $Q_1,\ldots, Q_s$ for some integer $s \ge 2$ and $Q_t$ be the last nontrivial component of $D$.
Then \[t < s\] by our assumption, $D_{1 \sim t}$ has the last strong component which is nontrivial, and $D_{t+1 \sim s}$ has only trivial strong components.
As we mentioned in the introduction of this section, $\left\{C^m(D_{t+1 \sim s})\right\}_{m=1}^{\infty}$ converges to an empty graph.
We note that all the arcs goes from a strong component of lower index to that of higher index.
Then the following is an immediate consequence of the above observations.

\begin{Prop}\label{prop:empty}
  Let $D$ be a multipartite tournament with ordered strong components $Q_1,\ldots, Q_s$ for some integer $s \ge 2$ and $Q_t$ ($t<s$) be the last nontrivial component of $D$.
  Then there is a positive integer $N$ such that the vertices in $D_{t+1 \sim s}$ are isolated in $C^m(D)$ for any integer $m \ge N$.
\end{Prop}

The following proposition is rather self-evident.

\begin{Prop}\label{prop:sub}
 Let $D'$ be a subdigraph of a digraph $D$.
 If $u$ and $v$ are adjacent in $C^m(D')$ for a positive integer $m$, then $u$ and $v$ are also adjacent in $C^m(D)$.
\end{Prop}

By Theorem~\ref{thm:last}, the graph sequence $\left\{C^m(D_{1 \sim t})\right\}_{m=1}^{\infty}$ converges and its limit is isomorphic to $K[D_{1 \sim t}]$ or one of the graphs given in Figure~\ref{fig:graphs}.
Accordingly, Proposition~\ref{prop:sub} tells us that,
\statement{onlyy}
{
	to investigate the graph sequence $\left\{C^m(D)\right\}_{m=1}^{\infty}$, we only need to take a look at nonadjacent vertices in the graphs in Figure~\ref{fig:graphs} corresponding to the limit of $\left\{C^m(D_{1 \sim t})\right\}_{m=1}^{\infty}$.
}
Furthermore, the following is true:
\statement{only}
{
	no two of nonadjacent vertices in the graphs corresponding to the limit of $\left\{C^m(D_{1 \sim t})\right\}_{m=1}^{\infty}$  in Figure~\ref{fig:graphs} have a step common prey in $D_{1 \sim t}$ (here, a step common prey of two vertices means an $m$-step common prey for some positive integer $m$).
}

By Proposition~\ref{prop:empty}, we may ignore the trivial strong components after the last nontrivial component when we study the graph sequence $\{C^m(D)\}_{m=1}^{\infty}$.
If the last nontrivial component is primitive,
then we can obtain the limit given in the following proposition by using Theorem~\ref{thm:last} and Propositions~\ref{prop:empty} and \ref{prop:sub}.

\begin{Prop}\label{prop:qt_primitive}
  Let $D$ be a multipartite tournament with ordered strong components $Q_1,\ldots, Q_s$ for some integer $s \ge 2$ and $Q_t$ ($t < s$) be the last nontrivial component of $D$.
  If $Q_t$ is primitive, then the graph sequence $\{C^m(D)\}_{m=1}^{\infty}$ converges to the graph isomorphic to  $K\left[D_{1 \sim t}\right] \cup I\left[D_{t+1 \sim s}\right]$.
\end{Prop}

\begin{Lem}\label{prop:notpartite-related}
  Let $D$ be a multipartite tournament with the ordered strong components $Q_1,\ldots, Q_s$ for some integer $s \ge 2$ and $Q_t$ ($t < s$) be the last nontrivial component of $D$.
  Suppose that there exists a trivial component $Q_j$ for some $j \in \{t+1, \ldots, s\}$ such that no partite sets of $Q_t$ are partite-related to the partite set of $Q_j$.
  Then the graph sequence $\{C^m(D)\}_{m=1}^{\infty}$ converges to the graph isomorphic to $K\left[D_{1 \sim t}\right] \cup I\left[D_{t+1 \sim s}\right]$.
\end{Lem}

\begin{proof}
  Take a vertex $v$ in $D_{1 \sim t}$ and fix a positive integer $m$.
  Since $Q_t$ is nontrivial, there are at least two partite sets in $Q_t$ and so there is an arc from $v$ to a vertex $x$ in $Q_t$.
  Then, since $Q_t$ is a nontrivial strong component, there is an $(x,y)$-directed walk of length $m-1$ for some vertex $y$ in $Q_t$ and therefore we obtain a $(v,y)$-directed walk of length $m$.
  By the hypothesis, $y \to w_j$ where $V(Q_j)= \{w_j\}$.
  Hence $w_j$ is an $(m+1)$-step prey of each vertex in $D_{1 \sim t}$ since $v$ was arbitrarily chosen.
  Since $m$ was chosen arbitrarily, we can conclude that the graph sequence $\{C^m(D)\}_{m=1}^{\infty}$ converges to the graph isomorphic to $K\left[D_{1 \sim t}\right] \cup I\left[D_{t+1 \sim s}\right]$ by Proposition~\ref{prop:empty}.
\end{proof}

Let  $\{C^m(D)\}_{m=1}^{\infty}$ be the graph sequence for a multipartite tournament $D$ with ordered strong components $Q_1,\ldots, Q_s$ for some integer $s \ge 2$ where $Q_t$ ($t < s$) is the last nontrivial component of $D$.
By \eqref{1234}, $\kappa(Q_t)\in \{1,2,3,4\}$.
By Proposition~\ref{prop:qt_primitive}, in order to compute the limit of the graph sequence $\{C^m(D)\}_{m=1}^{\infty}$, it remains to consider the multipartite tournaments with the last nontrivial component being imprimitive, that is, $\kappa(Q_t)= 2,3,4$.
We will take care of the cases $\kappa(Q_t)= 2$ and $\kappa(Q_t)=4$ first.

We need to introduce the notion of \textit{head completing index} of $D$, which is defined in the following way.
If there is an integer $r$ less than $t$ such that there is a partite set of $Q_r$ which is not partite-related to any partite set of $Q_t$, then we take the largest $r$ as the head completing index of $D$, that is, each partite set of $Q_{j}$ is partite-related to some partite set of $Q_t$ for each $j \in \{r+1, \ldots, t-1\}$.
Otherwise, we say that the head completing index of $D$ equals $0$.
For example, the head completing index of a bipartite tournament equals $0$.

In addition, for notational convenience, we let 
\begin{equation}\label{a12}
	 A_1=V_1 \cap \bigcup_{i=r+1}^{t}V(Q_i) \quad \mbox{and} \quad A_2=V_2 \cap \bigcup_{i=r+1}^{t}V(Q_i)
\end{equation}
where $r$ is the head completing index of $D$ and $V_1$ and $V_2$ are partite sets of $D$.

Now we consider a multipartite tournament $D$ with ordered strong components $Q_1,\ldots, Q_s$ for some integer $s \ge 2$ where $Q_t$ ($t<s$) is the last nontrivial component with  $\kappa(Q_t)=2$ or $4$.
We will compute the limit of the graph sequence $\{C^m(D)\}_{m=1}^{\infty}$.
To do so, it is necessary to first examine the limit when the last strong component is nontrivial, which is described in the proofs of Theorems~3.2 and 3.4 of Jung {\it{et al.}}~\cite{jung2023limit} as follows.

\begin{Prop}[\cite{jung2023limit}]\label{prop:head}
  Let $D$ be a $k$-partite tournament with $k$-partition $(V_1,\ldots,V_k)$ and ordered strong components $Q_1,\ldots, Q_s$ for some integers $k \ge 2$ and  $s \ge 2$ and let $Q_s$ be the last nontrivial component of $D$.
Then the graph sequence $\{C^m(D)\}_{m=1}^{\infty}$ converges to the graph isomorphic to
\begin{itemize}
  \item[(i)] the graph $G_1$ in Theorem~\ref{thm:last} if $\kappa(Q_s)=2$ where
    \[
     K^{(1)}=K\left[D_{1 \sim r}\right],
     K^{(2)}=K\left[A_1\right],
     \mbox{ and }
      K^{(3)}=K\left[A_2\right];
    \]
  \item[(ii)] the graph $G_2$ in the same theorem if $\kappa(Q_s)=4$ where
\[
K^{(1)} = K\left[D_{1 \sim r}\right],\
K^{(2)} = K\left[A_1 - V(Q_s) \right],\
K^{(3)} = K\left[A_2 - V(Q_s)\right],
\]
\[
K^{(4)} = K\left[U_{1}^{(s)}\right],\
K^{(5)} = K\left[U_{3}^{(s)}\right],\
K^{(6)} = K\left[U_{2}^{(s)}\right], \mbox{ and }
K^{(7)} = K\left[U_{4}^{(s)}\right] \]
\end{itemize}
where $r$ is the head completing index of $D$ and $A_1$ and $A_2$ are given in \eqref{a12}.
\end{Prop}

\begin{Thm}\label{thm:k=2}
Let $D$ be a $k$-partite tournament with $k$-partition $(V_1,\ldots,V_k)$ and ordered strong components $Q_1,\ldots, Q_s$ for some integers $k \ge 2$ and  $s \ge 2$ and let $Q_t$ be the last nontrivial component of $D$ with $\kappa(Q_t)=2$ for $t < s$.
Then the graph sequence $\{C^m(D)\}_{m=1}^{\infty}$ converges to the graph isomorphic to
\[  \begin{cases}
    \left(K\left[D_{1 \sim r}\right] \vee \left(K[A_1] \cup K[A_2] \right)\right)\cup I\left[D_{t+1 \sim s}\right], & \mbox{if $V(D_{t+1 \sim s}) \subseteq V_1 \cup V_2$};\\
    K\left[D_{1 \sim t}\right] \cup I\left[D_{t+1 \sim s}\right], & \mbox{otherwise}
  \end{cases}
\]
where $r$ is the head completing index of $D$ and $A_1$ and $A_2$ are given in \eqref{a12}.
\end{Thm}

\begin{proof}
	Suppose $V(D_{t+1 \sim s}) \not\subseteq V_1 \cup V_2$.
	Then there exists a trivial component $Q_i$ for some integer $i \in \{t+1,\ldots, s\}$ such that no partite sets of $Q_t$ are partite-related to the partite set of $Q_i$ by \eqref{eqn:bi}.
	Thus, by Lemma~\ref{prop:notpartite-related},
	the limit of the graph sequence $\{C^m(D)\}_{m=1}^{\infty}$ is isomorphic to $K\left[D_{1 \sim t}\right] \cup I\left[D_{t+1 \sim s}\right]$.
	
	Now suppose $V(D_{t+1 \sim s}) \subseteq V_1 \cup V_2$.
	Then $D_{r+1 \sim s}$ is a bipartite tournament.
	We note that, by Proposition~\ref{prop:empty}, there is a positive integer $N_1$ such that the vertices in $D_{t+1 \sim s}$ are isolated in $C^m(D)$ for any integer $m \ge N_1$.
	By Proposition~\ref{prop:head}(i), there exists a positive integer $N_2$ such that $C^m(D_{1 \sim t})$ is isomorphic to
	 $K\left[D_{1 \sim r}\right] \vee \left(K[A_1] \cup K[A_2] \right)$ for any integer $m \ge N_2$.

By \eqref{onlyy} and Proposition~\ref{prop:sub}, we only need to check the adjacency between a vertex in $A_1$ and a vertex in $A_2$.
Since $D_{r+1 \sim s}$ is a bipartite tournament, if there is a directed walk from a vertex in $A_1$ to a vertex in $V(D_{r+1 \sim s}) \cap V_1$ (resp.\ $V(D_{r+1 \sim s}) \cap V_2$), then its length is even (resp.\ odd) and
if there is a directed walk from a vertex in $A_2$ to a vertex in $V(D_{r+1 \sim s}) \cap V_1$ (resp.\ $V(D_{r+1 \sim s}) \cap V_2$), then its length is odd (resp.\ even).
Thus a vertex in $A_1$ and a vertex in $A_2$ do not have any $m$-step common prey in $D$ for any positive integer $m$ and so there are no edges joining a vertex in $A_1$ and a vertex in $A_2$ in $C^m(D)$.
Therefore $C^{m}(D)$ is isomorphic to $K\left[D_{1 \sim r}\right] \vee \left(K[A_1] \cup K[A_2] \right)\cup I\left[D_{t+1 \sim s}\right]$ for any integer $m \ge \max\{N_1, N_2\}$ and so the limit of the graph sequence $\{C^m(D)\}_{m=1}^{\infty}$ is isomorphic to $K\left[D_{1 \sim r}\right] \vee \left(K[A_1] \cup K[A_2] \right)\cup I\left[D_{t+1 \sim s}\right]$.
\end{proof}

\begin{Thm}\label{thm:4}
  Let $D$ be a $k$-partite tournament with $k$-partition $(V_1,\ldots,V_k)$ and ordered strong components $Q_1,\ldots, Q_s$ for some integers $k \ge 2$ and  $s \ge 2$ and let $Q_t$ be the last nontrivial component of $D$ with $\kappa(Q_t)=4$ for $t < s$.
  Then there exists a positive integer $N$ such that, for any $m \ge N$, $C^m(D)$ is isomorphic to
\begin{itemize}
  \item $K\left[D_{1 \sim t}\right] \cup I\left[D_{t+1 \sim s}\right]$,  if $V(D_{t+1 \sim s}) \not\subseteq V_1 \cup V_2$;

       \item $\left(K\left[D_{1 \sim r}\right] \vee \left(K[A_1] \cup K[A_2] \right) \right)\cup I\left[D_{t+1 \sim s}\right]$, if $V(D_{t+1 \sim s}) \subseteq V_1 \cup V_2$, $V(D_{t+1 \sim s}) \cap V_1 \neq \emptyset$, and $V(D_{t+1 \sim s}) \cap V_2 \neq \emptyset$;
  \item $G_4 \cup I\left[D_{t+1 \sim s}\right]$ for even $m$ and $G_5 \cup I\left[D_{t+1 \sim s}\right]$ for odd $m$, if $V(D_{t+1 \sim s}) \subseteq V_1$;
        \item $G_5 \cup I\left[D_{t+1 \sim s}\right]$ for even $m$ and $G_4 \cup I\left[D_{t+1 \sim s}\right]$ for odd $m$, if $V(D_{t+1 \sim s}) \subseteq V_2$;
\end{itemize}
where $r$ is the head completing index of $D$ and $A_1$ and $A_2$ are given in \eqref{a12} and $G_4$ and $G_5$ are the graphs given in Figure~\ref{fig:graphs_for4}.
 \begin{figure}
  \begin{center}
  \begin{tabular}{cc}
%
%
    \begin{tikzpicture}[auto,thick,scale=0.8]
    \tikzstyle{player}=[minimum size=5pt,inner sep=1pt,outer sep=0pt,circle,draw=black]
    \tikzstyle{player1}=[minimum size=5pt,inner sep=1pt,outer sep=3pt,rectangle,draw=white]
    \tikzstyle{source}=[minimum size=5pt,inner sep=0pt,outer sep=0pt,ball color=black, circle]
    \tikzstyle{arc}=[minimum size=5pt,inner sep=1pt,outer sep=1pt, font=\footnotesize]
    \path (90:3cm)     node [player1]  (a) {$K^{(1)}$};
    \path (300:3.6cm)       node [player1] (b){$K^{(3)}$};
    \path (240:3.6cm)       node [player1] (c){$K^{(2)}$};
    \path (0:2cm)       node [player1] (d){$K^{(6)}$};
    \path (0:4cm)       node [player1] (e){$K^{(7)}$};
    \path (180:2cm)       node [player1] (f){$K^{(5)}$};
    \path (180:4cm)       node [player1] (g){$K^{(4)}$};

    \draw[black,thick,-] (a) -- ++(b);
    \draw[black,thick,-] (a) -- ++(c);
    \draw[black,thick,-] (b) -- ++(d);
    \draw[black,thick,-] (b) -- ++(e);
    \draw[black,thick,-] (c) -- ++(f);
    \draw[black,thick,-] (c) -- ++(g);
    \draw[black,thick,-] (a) -- ++(f);
    \draw[black,thick,-] (a) -- ++(g);
    \draw[black,thick,-] (a) -- ++(d);
    \draw[black,thick,-] (a) -- ++(e);
\draw[black,thick,-] (g) -- ++(f);

    \end{tikzpicture}
    & \begin{tikzpicture}
\matrix [matrix of math nodes,left delimiter=(,right delimiter=),row sep=0.1cm,column sep=-0.1cm] (m) {
      J^{(1)} & J & J  & J & J  & J & J   \\
      J & J^{(2)} & O  & J & J  & O & O  \\
      J & O   &J^{(3)} & O & O  & J & J  \\
      J & J & O & J^{(4)} & J & O & O \\
      J & J & O &J   &J^{(5)}& O & O  \\
      J &  O & J  &  O   & O  & J^{(6)} & O  \\
      J & O & J & O & O & O & J^{(7)} \\
       };
\end{tikzpicture}\\
$G_4$ & $M_4$
\\
    &\\
    \\
    &\\

    \begin{tikzpicture}[auto,thick,scale=0.8]
    \tikzstyle{player}=[minimum size=5pt,inner sep=1pt,outer sep=0pt,circle,draw=black]
    \tikzstyle{player1}=[minimum size=5pt,inner sep=1pt,outer sep=3pt,rectangle,draw=white]
    \tikzstyle{source}=[minimum size=5pt,inner sep=0pt,outer sep=0pt,ball color=black, circle]
    \tikzstyle{arc}=[minimum size=5pt,inner sep=1pt,outer sep=1pt, font=\footnotesize]
    \path (90:3cm)     node [player1]  (a) {$K^{(1)}$};
    \path (300:3.6cm)       node [player1] (b){$K^{(3)}$};
    \path (240:3.6cm)       node [player1] (c){$K^{(2)}$};
    \path (0:2cm)       node [player1] (d){$K^{(6)}$};
    \path (0:4cm)       node [player1] (e){$K^{(7)}$};
    \path (180:2cm)       node [player1] (f){$K^{(5)}$};
    \path (180:4cm)       node [player1] (g){$K^{(4)}$};

    \draw[black,thick,-] (a) -- ++(b);
    \draw[black,thick,-] (a) -- ++(c);
    \draw[black,thick,-] (b) -- ++(d);
    \draw[black,thick,-] (b) -- ++(e);
    \draw[black,thick,-] (c) -- ++(f);
    \draw[black,thick,-] (c) -- ++(g);
    \draw[black,thick,-] (a) -- ++(f);
    \draw[black,thick,-] (a) -- ++(g);
    \draw[black,thick,-] (a) -- ++(d);
    \draw[black,thick,-] (a) -- ++(e);
    \draw[black,thick,-] (d) -- ++(e);

    \end{tikzpicture}
    & \begin{tikzpicture}
\matrix [matrix of math nodes,left delimiter=(,right delimiter=),row sep=0.1cm,column sep=-0.1cm] (m) {
      J^{(1)} & J & J  & J & J  & J & J  \\
      J & J^{(2)} & O  & J & J  & O & O  \\
      J & O   &J^{(3)} & O & O  & J & J  \\
      J & J & O & J^{(4)} & O & O & O \\
      J & J & O &O   &J^{(5)}& O & O  \\
      J &  O & J  &  O   & O  & J^{(6)} & J  \\
      J & O & J & O & O & J & J^{(7)}\\
     };
\end{tikzpicture}\\
$G_5$ & $M_5$
\end{tabular}
\end{center}
\caption{The graphs $G_4$ and $G_5$ in Theorem~\ref{thm:4} where
$K^{(1)} = K\left[D_{1 \sim r}\right]$,
$K^{(2)} = K\left[A_1 - V(Q_t) \right]$,
$K^{(3)} = K\left[A_2 - V(Q_t) \right]$,
$K^{(4)} = K\left[U_{1}^{(t)}\right]$,
$K^{(5)} = K\left[U_{3}^{(t)}\right]$,
$K^{(6)} = K\left[U_{2}^{(t)}\right]$, and
$K^{(7)} = K\left[U_{4}^{(t)}\right]$
(refer to \eqref{figure} for other notations).
Here, $r$ is the head completing index of $D$ and $A_1$ and $A_2$ are given in \eqref{a12}.
}\label{fig:graphs_for4}
\end{figure}
\end{Thm}
\begin{proof}
	Suppose $V(D_{t+1 \sim s}) \not\subseteq V_1 \cup V_2$.
	Then, since $\kappa(Q_t) = 4$, no partite sets of $Q_t$ are partite-related to the partite set of $Q_j$ for any $j \in \{t+1, \ldots, s\}$.
	Thus, by Lemma~\ref{prop:notpartite-related}, there exists a positive integer $N$ such that $C^m(D)$ is isomorphic to $K\left[D_{1 \sim t}\right] \cup I\left[D_{t+1 \sim s}\right]$ for any $m \ge N$.
	
	Suppose $V(D_{t+1 \sim s}) \subseteq V_1 \cup V_2$.
	We note that, by Proposition~\ref{prop:empty}, there is a positive integer $N_1$ such that the vertices in $D_{t+1 \sim s}$ are isolated in $C^m(D)$ for any integer $m \ge N_1$.
		By Proposition~\ref{prop:head}(ii), there exists a positive integer $N_2$ such that $C^m(D_{1 \sim t})$ is isomorphic to a graph $G_2$ given in Theorem~\ref{thm:last} for any integer $m \ge N_2$ where
	\[
	K^{(1)} = K\left[D_{1 \sim r}\right],\
	K^{(2)} = K\left[A_1 - V(Q_t) \right],\
	K^{(3)} = K\left[A_2 - V(Q_t) \right],
	\]
	\[
	K^{(4)} = K\left[U_{1}^{(t)}\right],\
	K^{(5)} = K\left[U_{3}^{(t)}\right],\
	K^{(6)} = K\left[U_{2}^{(t)}\right], \mbox{ and }
	K^{(7)} = K\left[U_{4}^{(t)}\right]. \]
Now we fix an integer $m \ge \max\{N_1, N_2\}$.
By \eqref{onlyy} and Proposition~\ref{prop:sub}, we only need to check the adjacency between a vertex in $K^{(4)} = K\left[U_{1}^{(t)}\right]$ and a vertex in $K^{(5)} = K\left[U_{3}^{(t)}\right]$; a vertex in $K^{(6)} = K\left[U_{2}^{(t)}\right]$ and a vertex in $K^{(7)} = K\left[U_{4}^{(t)}\right]$; a vertex in $A_1$ and a vertex in $A_2$
 (note that $A_1 = V\left(K^{(2)}\right) \cup V\left(K^{(4)}\right) \cup  V\left(K^{(5)}\right)$, $A_2 =  V\left(K^{(3)}\right) \cup V\left(K^{(6)}\right) \cup V\left(K^{(7)}\right)$, and $V(A_1 \cup A_2) = V(D_{r+1 \sim t})$).

Since $V(D_{r+1 \sim s}) \subseteq V_1 \cup V_2$, 
$D_{r+1 \sim s}$ is a bipartite tournament.
Thus, if there is a directed walk from a vertex in $A_1$ to a vertex in $V(D_{r+1 \sim s}) \cap V_1$ (resp.\ $V(D_{r+1 \sim s}) \cap V_2$), then its length is even (resp.\ odd) and
if there is a directed walk from a vertex in $A_2$ to a vertex in $V(D_{r+1 \sim s}) \cap V_1$ (resp.\ $V(D_{r+1 \sim s}) \cap V_2$), then its length is odd (resp.\ even).
Thus there are no step common prey of a vertex in $A_1$ and a vertex in $A_2$ and so there are no edges joining a vertex in $A_1$ and a vertex in $A_2$ in $C^m(D)$.

Fix $i \in \{1,2\}$.
Take two vertices $u_i \in U_{i}^{(t)}$ and $v_i \in U_{i+2}^{(t)}$.
We note that $U_1^{(t)}$, $U_2^{(t)}$, $U_3^{(t)}$, $U_4^{(t)}$ are the sets of imprimitivity of $Q_t$.
Then the following are true:
\statement{v1}{
if $i$ and $m$ have different parities, then each of $u_i$ and $v_i$ has an $(m-1)$-step prey belonging to $U_{2}^{(t)} \cup U_{4}^{(t)}$;}
\statement{v2}{
if $i$ and $m$ have the same parities, then each of $u_i$ and $v_i$ has an $(m-1)$-step prey belonging to $U_{1}^{(t)} \cup U_{3}^{(t)}$.}
Further, for each $j=1,2$, if $V(D_{t+1 \sim s}) \cap V_j \neq \emptyset$, then
\statement{common}{
\centering{
$U_{3-j}^{(t)} \cup U_{5-j}^{(t)} \to V(D_{t+1 \sim s}) \cap V_{j}$.
}}
If $V(D_{t+1 \sim s}) \subseteq V_j$ for some $j \in \{1,2\}$,
then $V(D_{t+1 \sim s}) \cap V_{3-j} = \emptyset$.
Therefore, for any $j = 1,2$,
\statement{not}{
 if $V(D_{t+1 \sim s}) \subseteq V_j$, then there are no arcs from any vertex in $U_{j}^{(t)} \cup U_{j+2}^{(t)}$, which is a subset of $V_j$, to any vertex in $D_{t+1 \sim s}$.}
If $i$ and $m$ have different parities and $V(D_{t+1 \sim s}) \cap V_1 \neq \emptyset$, then each vertex in $V(D_{t+1 \sim s}) \cap V_1$ is an $m$-step common prey of $u_i$ and $v_i$ by \eqref{v1} and \eqref{common}.
If $i$ and $m$ have the same parities and $V(D_{t+1 \sim s}) \cap V_2 \neq \emptyset$, then
each vertex in $V(D_{t+1 \sim s}) \cap V_2$ is an $m$-step common prey of $u_i$ and $v_i$ by \eqref{v2} and \eqref{common}.

Based on this observation, we may conclude that
if $V(D_{t+1 \sim s}) \cap V_1 \neq \emptyset$ and $V(D_{t+1 \sim s}) \cap V_2 \neq \emptyset$, then $u_i$ and $v_i$ have an $m$-step common prey in $D$ and so $C^m(D)$ is isomorphic to $\left(K\left[D_{1 \sim r}\right] \vee \left(K[A_1] \cup K[A_2] \right) \right)\cup I\left[D_{t+1 \sim s}\right]$.

Consider the case $V(D_{t+1 \sim s}) \subseteq V_1$.
Then,  $V(D_{t+1 \sim s}) \cap V_1 \neq \emptyset$ and $V(D_{t+1 \sim s}) \cap V_2 = \emptyset$.
Thus, by the above observation, $u_i$ and $v_i$ have an $m$-step common prey in $D$ if $i$ and $m$ have different parities.
If $i$ and $m$ have the same parities, then, by \eqref{only}, \eqref{v2}, and \eqref{not}, $u_i$ and $v_i$ do not have any $a$-step common prey in $D$ for any positive integer $a$.
Hence $C^m(D)$ is isomorphic to $G_4$ for even $m$; $G_5$ for odd $m$.

In the case $V(D_{t+1 \sim s}) \subseteq V_2$, we may apply the same argument as for the previous case by interchanging $V_1$ and $V_2$ to conclude that $C^m(D)$ is isomorphic to $G_5$ for even $m$; $G_4$ for odd $m$.
\end{proof}

Now we consider the case $\kappa(Q_t)=3$.

Given a digraph $D$ and subdigraphs $D_1$ and $D_2$ of $D$, 
we write $D_1 \stackrel{\ell}{\rightsquigarrow} D_2$
if there is a directed walk of length $\ell$ in $D$ from every vertex in $D_1$ to every vertex in $D_2$ for some positive integer $\ell$.

\begin{Lem}\label{Lem:t-1}
  Let $D$ be a $k$-partite tournament with $k$-partition $(V_1,\ldots,V_k)$ and ordered strong components $Q_1,\ldots, Q_s$ for some integers $k \ge 3$ and  $s \ge 2$ and let $Q_t$ be the last nontrivial component of $D$ with $\kappa(Q_t)=3$ for $t < s$.
  Then $D_{1 \sim t-1}  \stackrel{\ell}{\rightsquigarrow} D_{t+1 \sim s}$ for any integer $\ell \ge 2$.
\end{Lem}
\begin{proof}
Fix $\ell \ge 2$ and take a vertex $u$ in $D_{1 \sim t-1}$ and a vertex $v$ in $D_{t+1 \sim s}$.
Suppose that $v \not\in V_1 \cup V_2 \cup V_3$.
Then $Q_t \to v$.
Obviously, $u \not\in V_a$ for some $a \in \{1,2,3\}$.
Since $\kappa(Q_t)=3$, $Q_t \cap V_a \neq \emptyset$ by \eqref{eqn:bi} and so $u \to Q_t \cap V_a$.
Since $Q_t$ is nontrivial and strongly connected, $u$ can reach some vertex in $Q_{t}$ by traversing a directed walk of length $\ell-1$.
By the fact $Q_t \to v$, $v$ is an $\ell$-step prey of $u$.

Now suppose $v \in V_1 \cup V_2 \cup V_3$.
In the case $v \in V_{\alpha}$ for some $\alpha \in \{1,2,3\}$, we consider the following three cases:
$u \in V_{\alpha+1}$; $u \in V_{\alpha+2}$; $u \not\in V_{\alpha+1} \cup V_{\alpha+2}$ where all the subscripts are reduced to modulo $3$.
In each of the above three cases corresponding to $v \in V_{\alpha}$, there is a directed walk from $u$ to $v$ of each of lengths $3a-1$, $3a$, and $3a+1$ for each positive integer $a$, which is justified by Table $\alpha$ for each $\alpha = 1,2,3$.
The ordered pair $(i,j)$ in the table means that a desired directed walk can be obtained by $u$ entering $U_i^{(t)}$, rotating the sets of imprimitivity of $Q_t$ for $a-1$ revolutions and then moving further on to $U_j^{(t)}$ to leave it.
Therefore $v$ is an $\ell$-step prey of $u$.

\begin{center}
  \begin{tabular}{c|ccc}

  & $3a-1$ & \phantom{d}$3a$\phantom{d} & $3a+1$\\
\hline
  $u \in V_2$                 & $(3,3)$ & $(1,2)$ & $(3,2)$\\
  $u \in V_3$                & $(2,2)$ & $(2,3)$   & $(1,3)$ \\
 $u \not\in V_2 \cup V_3$      & $(2,2)$ & $(2,3)$ & $(3,2)$\\
 \multicolumn{4}{c}{}\\[-0.5em]
 \multicolumn{4}{c}{Table 1: $v \in V_1$}\\
 \end{tabular}
 
 \begin{tabular}{c|ccc}

  & $3a-1$ & \phantom{d}$3a$\phantom{d} & $3a+1$\\
\hline
  $u \in V_1$                 & $(3,3)$ & $(2,3)$ & $(2,1)$\\
  $u \in V_3$                & $(1,1)$ & $(2,3)$   & $(2,1)$ \\
 $u \not\in V_1 \cup V_3$      & $(1,1)$ & $(3,1)$ & $(1,3)$\\
   \multicolumn{4}{c}{}\\[-0.5em]
  \multicolumn{4}{c}{Table 2: $v \in V_2$}\\
 \end{tabular}
 
 \begin{tabular}{c|ccc}

  & $3a-1$ & \phantom{d}$3a$\phantom{d} & $3a+1$\\
\hline
  $u \in V_1$                 & $(2,2)$ & $(3,1)$ & $(3,2)$\\
  $u \in V_2$                & $(1,1)$ & $(1,2)$   & $(3,2)$ \\
 $u \not\in V_1 \cup V_2$      & $(1,1)$ & $(1,2)$ & $(2,1)$\\
   \multicolumn{4}{c}{}\\[-0.5em]
  \multicolumn{4}{c}{Table 3: $v \in V_3$}\\
 \end{tabular}
\end{center}

Since $u$, $v$, and $\ell$ were arbitrarily chosen, we have shown that every vertex in $D_{t+1 \sim s}$ is an $\ell$-step prey of each vertex in $D_{1 \sim t-1}$ for any integer $\ell \ge 2$.
\end{proof}

Now we are ready to take care of the case $\kappa(Q_t)=3$ and eventually prove Theorem~\ref{Thm:diverse}.

\begin{Thm}\label{thm:3}
  Let $D$ be a $k$-partite tournament with $k$-partition $(V_1,\ldots,V_k)$ and ordered strong components $Q_1,\ldots, Q_s$ for some integers $k \ge 3$ and  $s \ge 2$ and let $Q_t$ be the last nontrivial component of $D$ with $\kappa(Q_t)=3$ for $t < s$.
  Then there exists a positive integer $N$ such that, for any $m \ge N$ and any  $i \in \{1,2,3\}$ satisfying $i \equiv m-j+1 \pmod 3$,  $C^m(D)$ is isomorphic to

	\begin{itemize}
   \item $G_{1,i}\cup I\left[D_{t+1 \sim s}\right]$, if $V(D_{t+1 \sim s}) \subseteq V_j$ for some $j \in \{1,2,3\}$;
  \item $G_{2,i}\cup I\left[D_{t+1 \sim s}\right]$, if $D_{t+1 \sim s}$ satisfies the following property;
  \statement{dag}
  {there is an integer $\alpha \in \{t+1, \ldots, s-1\}$ such that $V(D_{t+1 \sim \alpha })$ and $V(D_{\alpha+1 \sim s})$ are partite-related to $U_{j}^{(t)}$ and $U_{j+1}^{(t)}$, respectively, for some $j \in \{1,2,3\}$ (identifying $U_{4}^{(t)}$ with $U_{1}^{(t)}$)}
   \item $K\left[D_{1 \sim t}\right] \cup I\left[D_{t+1 \sim s}\right]$, otherwise
\end{itemize}
where $G_{1,i}$ are the graphs given in Figure~\ref{fig:graphs_for3},
$V(G_{2,i}) = V(G_{1,i})$, 
$E(G_{2,1})=E(G_{1,1}) \cup E \left(K^{(5)}\vee (K^{(4)} \cup K^{(7)})\right)$,
$E(G_{2,2}) = E(G_{1,2})\cup E\left(K^{(7)}\vee (K^{(2)} \cup K^{(6)})\right)$, and
$E(G_{2,3})=E(G_{1,3})\cup E\left(K^{(6)}\vee (K^{(3)} \cup K^{(5)})\right)$.
\end{Thm}

\begin{figure}
  \begin{center}
  \begin{tabular}{ccc}
\begin{tikzpicture}[auto,thick,scale=0.75]
    \tikzstyle{player}=[minimum size=5pt,inner sep=1pt,outer sep=0pt,circle,draw=black]
    \tikzstyle{player1}=[minimum size=5pt,inner sep=0pt,outer sep=3pt,rectangle,draw=white]
    \tikzstyle{source}=[minimum size=5pt,inner sep=0pt,outer sep=0pt,ball color=black, circle]
    \tikzstyle{arc}=[minimum size=5pt,inner sep=1pt,outer sep=1pt, font=\footnotesize]
    \path (40:5cm)  node [player1]  (a) {$K^{(4)}$};
    \path (90:1.5cm)     node [player1]  (b) {$K^{(3)}$};
    \path (140:5cm)  node [player1]  (c) {$K^{(2)}$};
    \path (340:3.6cm)   node [player1]  (d) {$K^{(7)}$};
    \path (270:2cm)       node [player1]  (e){$K^{(6)}$};
    \path (200:3.6cm)  node [player1]  (f){$K^{(5)}$};
    \path (90:5cm)  node [player1]  (g){$K^{(1)}$};

    \draw[black,thick,-] (a) -- ++(b);
    \draw[black,thick,-] (b) -- ++(c);
   \draw[black,thick,-] (a) -- ++(c);
    \draw[black,thick,-] (a) -- ++(e);
    \draw[black,thick,-] (c) -- ++(e);
    \draw[black,thick,-] (b) -- ++(f);
    \draw[black,thick,-] (b) -- ++(d);
    \draw[black,thick,-] (c) -- ++(f);
    \draw[black,thick,-] (a) -- ++(d);
    \draw[black,thick,-] (g) -- ++(a);
    \draw[black,thick,-] (g) -- ++(b);
    \draw[black,thick,-] (g) -- ++(c);
    \draw[black,thick,-] (g) -- ++(d);
   \path (g) edge [black, bend right=23,thick,-] (e);
    \draw[black,thick,-] (g) -- ++(f);
    \draw[black,thick,-] (d) -- ++(e);
    \draw[black,thick,-] (b) -- ++(e);
     \path (c) edge [black, bend right=16,thick,-] (d);
    \end{tikzpicture}
         & &
        \begin{tikzpicture}[scale=0.95, every node/.style={scale=0.95}]
\matrix [matrix of math nodes,left delimiter=(,right delimiter=),row sep=0.1cm,column sep=-0.1cm] (m) {
      J^{(1)} & J &   J  & J & J   & J &J  \\
      J & J^{(2)} &   J  & J & J   & J & J  \\
      J & J   &J^{(3)} & J & J   & J & J  \\
      J & J & J & J^{(4)} & O & J & J  \\
      J & J    & J    &   O   &J^{(5)}& O & O  \\
      J &  J    & J    &  J   & O  & J^{(6)} &J  \\
      J & J & J & J & O & J & J^{(7)} \\
     };
\end{tikzpicture}\\
        $G_{1,1}$  & & $A(G_{1,1})+I$
        \\\\
        \begin{tikzpicture}[auto,thick,scale=0.75]
    \tikzstyle{player}=[minimum size=5pt,inner sep=1pt,outer sep=0pt,circle,draw=white]
    \tikzstyle{player1}=[minimum size=5pt,inner sep=0pt,outer sep=3pt,rectangle,draw=white]
    \tikzstyle{source}=[minimum size=5pt,inner sep=0pt,outer sep=0pt,ball color=black, circle]
    \tikzstyle{arc}=[minimum size=5pt,inner sep=1pt,outer sep=1pt, font=\footnotesize]
    \path (40:5cm)  node [player1]  (a) {$K^{(4)}$};
    \path (90:1.5cm)     node [player1]  (b) {$K^{(3)}$};
    \path (140:5cm)  node [player1]  (c) {$K^{(2)}$};
    \path (340:3.6cm)   node [player1]  (d) {$K^{(7)}$};
    \path (270:2cm)       node [player1]  (e){$K^{(6)}$};
    \path (200:3.6cm)  node [player1]  (f){$K^{(5)}$};
    \path (90:5cm)  node [player1]  (g){$K^{(1)}$};

    \draw[black,thick,-] (a) -- ++(b);
    \draw[black,thick,-] (b) -- ++(c);
   \draw[black,thick,-] (a) -- ++(c);
    \draw[black,thick,-] (a) -- ++(e);
    \draw[black,thick,-] (c) -- ++(e);
    \draw[black,thick,-] (b) -- ++(f);
    \draw[black,thick,-] (b) -- ++(d);
    \draw[black,thick,-] (c) -- ++(f);
    \draw[black,thick,-] (a) -- ++(d);
    \draw[black,thick,-] (g) -- ++(a);
    \draw[black,thick,-] (g) -- ++(b);
    \draw[black,thick,-] (g) -- ++(c);
    \draw[black,thick,-] (g) -- ++(d);
   \path (g) edge [black, bend right=23,thick,-] (e);
    \draw[black,thick,-] (g) -- ++(f);
    \draw[black,thick,-] (e) -- ++(f);
    \draw[black,thick,-] (b) -- ++(e);
     \path (a) edge [black, bend left=16,thick,-] (f);
    \end{tikzpicture}
         & &
        \begin{tikzpicture}[scale=0.95, every node/.style={scale=0.95}]
\matrix [matrix of math nodes,left delimiter=(,right delimiter=),row sep=0.1cm,column sep=-0.1cm] (m) {
      J^{(1)} & J &   J  & J & J   & J &J \\
      J & J^{(2)} &   J  & J & J   & J & O  \\
      J & J   &J^{(3)} & J & J   & J & J  \\
      J & J & J & J^{(4)} & J & J & J \\
      J & J    & J    &   J   &J^{(5)}& J & O  \\
      J &  J    & J    &  J   & J  & J^{(6)} &O  \\
      J & O & J & J & O & O & J^{(7)} \\
     };
\end{tikzpicture}\\
        $G_{1,2}$  & & $A(G_{1,2})+I$
        \\\\
        \begin{tikzpicture}[auto,thick,scale=0.75]
    \tikzstyle{player}=[minimum size=5pt,inner sep=1pt,outer sep=0pt,circle,draw=white]
    \tikzstyle{player1}=[minimum size=5pt,inner sep=0pt,outer sep=3pt,rectangle,draw=white]
    \tikzstyle{source}=[minimum size=5pt,inner sep=0pt,outer sep=0pt,ball color=black, circle]
    \tikzstyle{arc}=[minimum size=5pt,inner sep=1pt,outer sep=1pt, font=\footnotesize]
   \path (40:5cm)  node [player1]  (a) {$K^{(4)}$};
    \path (90:1.5cm)     node [player1]  (b) {$K^{(3)}$};
    \path (140:5cm)  node [player1]  (c) {$K^{(2)}$};
    \path (340:3.6cm)   node [player1]  (d) {$K^{(7)}$};
    \path (270:2cm)       node [player1]  (e){$K^{(6)}$};
    \path (200:3.6cm)  node [player1]  (f){$K^{(5)}$};
    \path (90:5cm)  node [player1]  (g){$K^{(1)}$};

    \draw[black,thick,-] (a) -- ++(b);
    \draw[black,thick,-] (b) -- ++(c);
   \draw[black,thick,-] (a) -- ++(c);
    \draw[black,thick,-] (a) -- ++(e);
    \draw[black,thick,-] (c) -- ++(e);
    \draw[black,thick,-] (b) -- ++(f);
    \draw[black,thick,-] (b) -- ++(d);
    \draw[black,thick,-] (c) -- ++(f);
    \draw[black,thick,-] (a) -- ++(d);
    \draw[black,thick,-] (g) -- ++(a);
    \draw[black,thick,-] (g) -- ++(b);
    \draw[black,thick,-] (g) -- ++(c);
    \draw[black,thick,-] (g) -- ++(d);
   \path (g) edge [black, bend right=23,thick,-] (e);
    \draw[black,thick,-] (g) -- ++(f);
    \draw[black,thick,-] (d) -- ++(f);
     \path (a) edge [black, bend left=16,thick,-] (f);
     \path (c) edge [black, bend right=16,thick,-] (d);
    \end{tikzpicture}
         & &
        \begin{tikzpicture}[scale=0.95, every node/.style={scale=0.95}]
\matrix [matrix of math nodes,left delimiter=(,right delimiter=),row sep=0.1cm,column sep=-0.1cm] (m) {
      J^{(1)} & J &   J  & J & J   & J &J \\
      J & J^{(2)} &   J  & J & J   & J & J \\
      J & J   &J^{(3)} & J & J   & O & J  \\
      J & J & J & J^{(4)} & J & J & J \\
      J & J    & J    &   J   &J^{(5)}& O & J  \\
      J &  J    & O    &  J   & O  & J^{(6)} &O  \\
      J & J & J & J & J & O & J^{(7)}\\
       };
\end{tikzpicture}\\
        $G_{1,3}$ & &  $A(G_{1,3})+I$
\end{tabular}
\end{center}
\caption{The graphs $G_{1,1}$, $G_{1,2}$, and $G_{1,3}$ in Theorem~\ref{thm:3} (refer to \eqref{figure} for other notations).
}\label{fig:graphs_for3}
\end{figure}

\begin{proof}
	We take a positive integer $N$ such that $C^m(D_{1 \sim t})$ is isomorphic to a graph $G_3$ given in Theorem~\ref{thm:last} for any integer $m \ge N$ and satisfies the property described in Proposition~\ref{prop:empty}.
	We may assume $N \ge 4$.
	Now we fix an integer $m \ge N$.
It is sufficient to examine the adjacency between the vertices in $D_{1 \sim t}$.

 Suppose $V(D_{t+1 \sim s}) \not\subseteq V_1 \cup V_2 \cup V_3$.
 Then, since $\kappa(Q_t) = 3$, no partite sets of $Q_t$ are partite-related to the partite set of $Q_j$ for any $j \in \{t+1, \ldots, s\}$.
 Thus $C^m(D)$ is isomorphic to $K\left[D_{1 \sim t}\right] \cup I\left[D_{t+1 \sim s}\right]$ by Lemma~\ref{prop:notpartite-related}.
 
Suppose \[V(D_{t+1 \sim s}) \subseteq V_1 \cup V_2 \cup V_3.\]
Let \[\nu=\left|\{i \mid V(D_{t+1 \sim s}) \cap V_i \neq \emptyset\}\right|.\]
Then $\nu \in \{1,2,3\}$.
We first handle the case $\nu = 3$, which is rather simple to take care of.
Then $V(D_{t+1 \sim s}) \cap V_i \neq \emptyset$ for each $i=1,2,3$.
Take $u, v$ in $D_{1 \sim t}$.
Since $\kappa(Q_t)=3$, $u$ and $v$ can reach some vertices in $Q_t$, say $x$ and $y$, respectively, in $D_{1 \sim t}$ by traversing directed walks of length $m-1$.
Then, for some $a \in \{1,2,3\}$, none of $x$ and $y$ belongs to $V_a$.
By the case assumption, $V(D_{t+1 \sim s}) \cap V_a \neq \emptyset$.
Thus $x \to V(D_{t+1 \sim s}) \cap V_a$ and $y \to V(D_{t+1 \sim s}) \cap V_a$ and so $u$ and $v$ have an $m$-step common prey in $V(D_{t+1 \sim s}) \cap V_a$.
Since $u$ and $v$ were arbitrarily chosen, the subgraph of $C^m(D)$ induced by $V(D_{1 \sim t})$ is isomorphic to $K\left[D_{1 \sim t}\right]$.
Hence $C^m(D)$ is isomorphic to $K\left[D_{1 \sim t}\right] \cup I\left[D_{t+1 \sim s}\right]$.

In the following, we take care of the cases $\nu = 1$ or $2$.
By~\eqref{onlyy},
\statement{onlycheck}{
we only need to check the adjacency among $K^{(5)}$, $K^{(6)}$, and $K^{(7)}$; between $K^{(2)}$ and $K^{(7)}$; between $K^{(3)}$ and $K^{(6)}$; between $K^{(4)}$ and $K^{(5)}$ in the graph $G_3$ in Figure~\ref{fig:graphs}.
}
We note that $G_3$ is the same as $G_3$ in  \cite{jung2023limit} in which it is shown that
\begin{equation}\label{eq:234}
	V(K^{(1)} \cup K^{(2)} \cup K^{(3)} \cup K^{(4)}) \subseteq V(D_{1 \sim t-1})
\end{equation}
 and
\begin{equation}\label{eq:567}
K^{(5)} = K\left[U_{1}^{(t)}\right], \quad K^{(6)} = K\left[U_{2}^{(t)}\right], \quad \mbox{ and } \quad K^{(7)} = K\left[U_{3}^{(t)}\right]
\end{equation}
(see the proofs of Theorems~4.4 and 4.5 in \cite{jung2023limit}).
For any pair of cliques mentioned in \eqref{onlycheck}, any vertex in one clique and any vertex in the other clique have potential step common prey in $D_{t+1 \sim s}$ by \eqref{only}.
By Lemma~\ref{Lem:t-1}, $D_{1 \sim t-1}  \stackrel{\ell}{\rightsquigarrow} D_{t+1 \sim s}$ for any integer $\ell \ge 2$.
We note that $V(Q_t) = V(K^{(5)}) \cup V(K^{(6)}) \cup V(K^{(7)})$ by \eqref{eq:567}.
Then, by \eqref{eq:234} and the observation in \eqref{onlycheck}, we just need to see if there is a directed walk of a certain length from a vertex in $Q_t$ to a vertex in $D_{t+1 \sim s}$.

{\it Case 1.}  $\nu = 1$, i.e.,\ $V(D_{t+1 \sim s}) \subseteq V_j$ for some $j \in \{1,2,3\}$.
Take a positive integer $a$.
Then, in the table given in \eqref{table:q_t1}, O (resp.\ X) in the $(p,q)$-entry indicates the existence (resp.\ non-existence) of a directed walk of the length at the top of the $q$th column from a vertex in the set at the leftmost of the $p$th row to a vertex in  $D_{t+1 \sim s}$ (given in the header row) for each $p,q=1,2,3$.
For example, the first row tells us that any vertex in $D_{t+1 \sim s}$ is a $(3a+j-2)$-step prey and a $(3a+j-1)$-step prey of any vertex in $U_{1}^{(t)}$ while no vertices in $D_{t+1 \sim s}$ are a $(3a+j)$-step prey of a vertex in $U_{1}^{(t)}$.
\begin{equation}\label{table:q_t1}
\begin{tabular}{c|ccc}
  & \multicolumn{3}{c}{$D_{t+1 \sim s}$} \\
  \cline{2-4}
   & $3a+j-2$ & \phantom{d}$3a+j-1$\phantom{d} & $3a+j$ \\
\hline
  $U_{1}^{(t)}$ & O & O & X \\
  $U_{2}^{(t)}$ & O & X & O \\
  $U_{3}^{(t)}$ & X & O & O \\
 \end{tabular}
\end{equation}

Let $\{K^{(p)}, K^{(q)} \}$ with $p < q$ 
be a pair taken from the cliques whose adjacency in $G_3$ given in Figure~\ref{fig:graphs} has to be checked based upon \eqref{onlycheck}.
Then \[(p,q) \in \{(2,7), (3,6), (4,5), (5,6), (5,7), (6,7)\}.\]
Based on \eqref{eq:567}, \eqref{table:q_t1} and Lemma~\ref{Lem:t-1}, we may conclude the following.
If a pair $(p,q)$ is contained in $\{ (3,6), (4,5), (5,6) \}$, then any vertex in $K^{(p)}$ and any vertex in $K^{(q)}$ are adjacent in $C^{3a+j-2}(D)$.
For any vertex in $K^{(7)}$, it cannot have a vertex in $D_{t+1 \sim s}$ as a $(3a+j-2)$-step prey and so cannot be adjacent to any vertex of $K^{(2)}$ or $K^{(5)}$ or $K^{(6)}$ in $C^{3a+j-2}(D)$.
Thus we have shown that $C^m(D)$ is isomorphic to $G_{1,2}$ given in Figure~\ref{fig:graphs_for3} if $m = 3a+j-2$.
Note that $m-j+1 \equiv (3a+j-2)-j+1 \equiv 2 \pmod{3}$ and so $i=2$.
Similarly, one may show that if $m = 3a+j$ (resp.\ $m = 3a+j-1$), then $C^{m}(D)$ is isomorphic to $G_{1,1}$ (resp.\ $G_{1,3}$) given in Figure~\ref{fig:graphs_for3}.

{\it Case 2.} $\nu = 2$ and $D_{t+1 \sim s}$ satisfies the property~\eqref{dag}, that is, there is an integer $\alpha \in \{t+1, \ldots, s-1\}$ such that $V(D_{t+1 \sim \alpha })$ and $V(D_{\alpha+1 \sim s})$ are partite-related to $U_{j}^{(t)}$ and $U_{j+1}^{(t)}$, respectively, for some $j \in \{1,2,3\}$ (identifying $U_{4}^{(t)}$ with $U_{1}^{(t)}$).
Take a positive integer $a$.
Now, to see the presence of a directed walk from a vertex in $Q_t$ to a vertex in $D_{t+1 \sim s}$ based on the length, we refer to the following tables.
We note that, except for the header row, the table given in  \eqref{table:q_t} is identical to the table given in \eqref{table:q_t1}.
Therefore $C^m(D)$ contains a subgraph isomorphic to $G_{1,i}$ where $i \equiv m-j+1 \pmod 3$ and $i \in \{1,2,3\}$.
\begin{equation}\label{table:q_t}
\begin{tabular}{c|ccc}
  & \multicolumn{3}{c}{$D_{t+1 \sim \alpha}$} \\
  \cline{2-4}
   & $3a+j-2$ & $\phantom{d}3a+j-1\phantom{d}$ & $3a+j$ \\
\hline
  $U_{1}^{(t)}$ & O & O & X \\
  $U_{2}^{(t)}$ & O & X & O \\
  $U_{3}^{(t)}$ & X & O & O \\
 \end{tabular}
\end{equation}
\begin{equation}\label{table:q_t2}
 \begin{tabular}{c|ccc}
  & \multicolumn{3}{c}{$D_{\alpha +1 \sim s}$} \\
  \cline{2-4}
   & $3a+j-2$ & $\phantom{d}3a+j-1\phantom{d}$ & $3a+j$ \\
\hline
  $U_{1}^{(t)}$ & X & O & O \\
  $U_{2}^{(t)}$ & O & O & X \\
  $U_{3}^{(t)}$ & O & X & O \\
 \end{tabular}
\end{equation}

For the table given in \eqref{table:q_t2}, O in the $(2,1)$-entry and O in the $(3,1)$-entry indicate that every vertex in $K^{(6)}$ and every vertex in $K^{(7)}$ have a common $(3a+j-2)$-step prey in $D_{\alpha + 1 \sim s}$.
Therefore every vertex in $K^{(6)}$ and every vertex in $K^{(7)}$ are adjacent in $C^{3a+j-2}(D)$.
Furthermore, every vertex in $K^{(2)}$ and every vertex in $K^{(7)}$ are adjacent in $C^{3a+j-2}(D)$ by Lemma~\ref{Lem:t-1}.
On the other hand, no vertices in $K^{(5)}$ and no vertices in $K^{(7)}$ have a $(3a+j-2)$-step common prey in $D_{t+1 \sim \alpha}$ or $D_{\alpha +1 \sim s}$ as evident from the tables in \eqref{table:q_t} and \eqref{table:q_t2}.
Therefore $C^{3a+j-2}(D)$ is isomorphic to $G_{2,2}$ (recall $E(G_{2,2})=E(G_{1,2})\cup E\left(K^{(7)}\vee (K^{(2)} \cup K^{(6)})\right)$).
For similar reasons, if $m = 3a+j$ (resp.\ $m = 3a+j-1$), then $C^{m}(D)$ is isomorphic to $G_{2,1}$ (resp.\ $G_{2,3}$).

{\it Case 3.} $\nu = 2$ and $D_{t+1 \sim s}$ does not satisfy the property~\eqref{dag}.
Since $D_{t+1 \sim s}$ does not satisfy the property~\eqref{dag}, there is an arc from a vertex $z \in V(D_{t+1 \sim s}) \cap V_{j+1}$ to a vertex $w \in V(D_{t+1 \sim s}) \cap V_j$ for some $j \in \{1,2,3\}$ (identifying $V_{4}^{(t)}$ and $V_{5}^{(t)}$ with $V_{1}^{(t)}$ and  $V_{2}^{(t)}$, respectively).

Starting at any arbitrary vertex of $U_{i}^{(t)}$ for each $i = 1,2,3$, it is possible to create a directed walk of any length greater than or equal to $2$ toward $w$ by traversing the necessary number of sets of imprimitivity of $Q_t$.
For example, we consider the case $i=j$ for $j$ given above.
To create a directed walk from a vertex $v$ in $U_{j}^{(t)}$ to $w$ of length $3a-1$ (resp.\ $3a$) for a positive integer $a$, 
\begin{itemize}
	\item we start from $v$;
	\item rotate the sets of imprimitivity of $Q_t$ for $a-1$ revolutions;
	\item move by one unit (resp.\ two units) further on to $U_{j+1}^{(t)} (\subseteq V_{j+1})$ (resp.\ $U_{j+2}^{(t)} (\subseteq V_{j+2})$) where $U_{4}^{(t)}$ and $U_{5}^{(t)}$ are identified with $U_{1}^{(t)}$ and  $U_{2}^{(t)}$, respectively;
	\item we can reach $w$ by passing through an arc (note that $w \in V_{j}$).
\end{itemize}
Therefore, there exists a directed walk from each vertex of $U_{j}^{(t)}$ to $w$ having a length of $3a-1$ or $3a$.
To create a directed walk from a vertex $v$ in $U_j^{(t)}$ to $w$ of length $3a+1$ for a positive integer $a$,
\begin{itemize}
	\item we start from $v$;
	\item rotate the sets of imprimitivity of $Q_t$ for $a-1$ revolutions;
	\item move by two units further on to $U_{j+2}^{(t)} (\subseteq V_{j+2})$  where $U_{4}^{(t)}$ and $U_{5}^{(t)}$ are identified with $U_{1}^{(t)}$ and  $U_{2}^{(t)}$, respectively;
	\item we can reach $z$ by passing through an arc (note that $z \in V_{j+1}$) and then $w$ by passing through an arc $z \to w$.
\end{itemize}
Since there is an arc from $z$ to $w$, there exists a directed walk from each vertex of $U_{j}^{(t)}$ to $w$ with a length of $3a+1$.
Similarly, one can show that there is a directed walk from a vertex in $U_{j+1}^{(t)}$ or a vertex in $U_{j+2}^{(t)}$ to $w$ of any length greater than or equal to 2.
By doing so, we can come to the conclusion that $Q_t  \stackrel{m}{\rightsquigarrow} w$ for any integer $m \ge 2$.

Furthermore, by Lemma~\ref{Lem:t-1}, $D_{1 \sim t-1}  \stackrel{\ell}{\rightsquigarrow} D_{t+1 \sim s}$ for any integer $\ell \ge 2$.
Therefore $w$ is an $m$-step common prey of the vertices in $D_{1 \sim t}$ and so $C^{m}(D)$ is isomorphic to $K\left[D_{1 \sim t}\right] \cup I\left[D_{t+1 \sim s}\right]$ for each integer $m \ge N$.
\end{proof}

\section{Acknowledgement}
This research was supported by the National Research Foundation of Korea(NRF)  (NRF-2022R1A2C1009648 and NRF-2017R1E1A1A03070489) funded by the Korea government(MSIP).

\bibliographystyle{plain}

\end{document}